\documentclass[11pt,a4paper]{article}

\usepackage{booktabs}
\usepackage{amsfonts}
\usepackage{amsmath,amsthm}
\usepackage{graphicx}
\usepackage{algorithmic}
\usepackage[ruled]{algorithm2e}
\SetAlFnt{\small}
\SetAlCapFnt{\small}
\SetAlCapNameFnt{\small}
\SetAlCapHSkip{0pt}
\IncMargin{-\parindent}
\usepackage{tikz,pgfplots}
\usepackage[hidelinks]{hyperref}
\usepackage[left=2.5cm,right=2.5cm,top=3cm,bottom=3cm]{geometry}
\usepackage{cleveref}
\Crefname{ALC@unique}{Line}{Lines}

\newcounter{myalg}
\makeatletter
\@addtoreset{ALC@unique}{myalg}
\makeatother

\newenvironment{acknowledgements}{\paragraph{Acknowledgements}}{}
\usepackage{amsopn}
\DeclareMathOperator{\diag}{diag}

\title{Thick-restarted joint Lanczos bidiagonalization for the GSVD\thanks{This work was supported by the Spanish Agencia Estatal de Investigaci{\'o}n under grant PID2019-107379RB-I00 / AEI / 10.13039/501100011033.}
}

\author{Fernando Alvarruiz\thanks{D.\ Sistemes Inform\`atics i Computaci\'o, Universitat Polit\`ecnica de Val\`encia, Val\`encia, Spain
  (\texttt{fbermejo@dsic.upv.es}).}
\and Carmen Campos\thanks{D.\ Did\'actica de la Matem\'atica, Universitat de Val\`encia, Val\`encia, Spain
  (\texttt{carmen.campos-gonzalez@uv.es}).}
\and Jose E. Roman\thanks{D.\ Sistemes Inform\`atics i Computaci\'o, Universitat Polit\`ecnica de Val\`encia, Val\`encia, Spain
  (\texttt{jroman@dsic.upv.es}).}
}

\begin{document}

\maketitle

\begin{abstract}
The computation of the partial generalized singular value decomposition (GSVD) of large-scale matrix pairs can be approached by means of iterative methods based on expanding subspaces, particularly Krylov subspaces. We consider the joint Lanczos bidiagonalization method, and analyze the feasibility of adapting the thick restart technique that is being used successfully in the context of other linear algebra problems. Numerical experiments illustrate the effectiveness of the proposed method. We also compare the new method with an alternative solution via equivalent eigenvalue problems, considering accuracy as well as computational performance. The analysis is done using a parallel implementation in the SLEPc library.
\end{abstract}

\section{Introduction}\label{sec:intro}

The generalized singular value decomposition (GSVD) of two matrices was introduced by Van Loan~\cite{Loan:1976:GSV}, with subsequent additional developments by Paige and Saunders~\cite{Paige:1981:TGS}. Given two real matrices $A$ and $B$ with the same number of columns, the GSVD of the pair $\{A,B\}$ is given by
\begin{equation}
\label{eq:gsvd}
U_A^TAG=C,\qquad U_B^TBG=S,
\end{equation}
where $U_A$, $U_B$ are orthogonal matrices, $G$ is a square nonsingular matrix, and $C$, $S$ are diagonal matrices in the simplest case. A more precise definition of the problem will be given in \cref{sec:svd}. This factorization, as an extension of the usual singular value decomposition (SVD), is finding its way in an increasing number of applications, for instance in the solution of constrained least squares problems~\cite{Golub:1996:MC}, discrete ill-posed problems via Tikhonov regularization~\cite{Kilmer:2007:PAG,Jia:2020:JBB}, as well as in many other contexts~\cite{Edelman:2020:GWA}.

The problem of computing the GSVD of a small dense matrix pair is well understood, and a robust implementation is available in LAPACK~\cite[\S 2.3.5.3]{Anderson:1999:LUG}. However, the case of large sparse matrix pairs is still the subject of active research. Normally, the computation of the large-scale GSVD is addressed by means of iterative methods based on expanding subspaces. This is the case of the Jacobi--Davidson method proposed by Hochstenbach~\cite{Hochstenbach:2009:JTM}. In this paper, we focus on Krylov methods, which still need to incorporate a great deal of the knowledge and innovation that has been successfully applied in similar linear algebra problems such as the SVD. One example of such is an effective restart mechanism, which is the main focus of this paper. The ultimate goal is to provide reliable software implementations of the methods that are readily available to users from the field of scientific computing or any other discipline that requires computing a GSVD. In our case, the developments presented in this paper are included in one of the solvers of SLEPc, the Scalable Library for Eigenvalue Problem Computations~\cite{Hernandez:2005:SSF}.

Zha~\cite{Zha:1996:CGS} was the first to propose a Lanczos method to be used for the GSVD. His method relies on a joint bidiagonalization of the two matrices, $A$ and $B$, in a similar way as Lanczos methods for the SVD have a single bidiagonalization in their foundation. The joint bidiagonalization in Zha's method results in two small-sized upper bidiagonal matrices. Later, Kilmer and coauthors~\cite{Kilmer:2007:PAG} proposed a variant in which one of the two bidiagonal matrices generated by the joint bidiagonalization is lower instead of upper. The latter has since then been the most popular approach, and we focus mainly on that variant for our developments.

As in any Lanczos method, a finite-precision arithmetic implementation suffers from various numerical pitfalls, such as loss of orthogonality in the generated Krylov basis vectors. Jia and Li~\cite{Jia:2023:JBM} have studied numerical error in the context of joint bidiagonalization for the GSVD, showing that loss of orthogonality can be prevented by simply enforcing semiorthogonality of Lanczos vectors, e.g., with a partial reorthogonalization technique~\cite{Jia:2021:JBP}, as is done in Lanczos methods for other linear algebra problems. Either in the case that a semiorthogonality scheme is pursued, or a full reorthogonalization approach is followed as we do to avoid loss of orthogonality, a consequence is that all Lanczos vectors must be kept throughout the computation, with the consequent increase in storage and computational costs.

Thick restart is an effective mechanism to keep the size of the Krylov basis bounded, that has been applied to Lanczos methods in many different contexts such as the symmetric eigenproblem~\cite{Wu:2000:TLM} or the SVD~\cite{Baglama:2005:AIR,Hernandez:2008:REP}. Compared to explicit restart, the thick restart scheme is much more effective because it compresses the currently available Krylov subspace into another Krylov subspace of smaller dimension (not a single vector), that retains the wanted spectral information while purging the components associated with unwanted eigenvalues (or singular values). The thick restart methodology can be summarized in two stages. In the first stage, one or more Lanczos recurrences are used to build one or more sets of Lanczos vectors, in a way that the small-sized problem resulting from the projection of the original problem retains the properties of the original problem (structure preservation). For instance, for symmetric-indefinite matrix pencils it is possible to employ a pseudo-Lanczos recurrence that results in a symmetric-indefinite projected problem~\cite{Campos:2016:RQM}. In the second stage, the built factorization is truncated to a smaller size decomposition, in a way that it is feasible to extend it again using the same Lanczos recurrences.

In this paper, we work out the details that are needed for thick-restarting the Lanczos recurrences associated with the joint bidiagonalization for the GSVD, so that the projected problem is a small-size GSVD, and this structure is preserved whenever the restart truncates the involved decompositions and they are extended again.

The rest of the paper is organized as follows. \Cref{sec:background} presents all the background material that is required for \cref{sec:trlanczos}, which presents the new developments related to thick restart for the GSVD. In \cref{sec:implem} we provide a few details of how the proposed method has been implemented in the SLEPc library. \Cref{sec:results} illustrates how the solver performs when applied to several test problems. Finally, we wrap up with some concluding remarks in \cref{sec:concl}. Throughout the paper, the presentation is done for real matrices, although the extension to the complex case is straightforward. In fact, our implemented solver supports both real and complex arithmetic.

\section{Background}\label{sec:background}

In this section, we review a number of concepts that are required for the developments of subsequent sections. Many of the concepts are also discussed for the case of the SVD, aiming at facilitating the understanding of the GSVD case, which is more involved.

\subsection{The SVD and the GSVD}\label{sec:svd}

Recall that the (standard) SVD of a matrix $A\in\mathbb{R}^{m\times n}$ is written as
\begin{equation}
\label{eq:svd}
A=U\Sigma V^T,
\end{equation}
where $U=[u_1,\ldots,u_m]\in\mathbb{R}^{m\times m}$ and $V=[v_1,\ldots,v_n]\in\mathbb{R}^{n\times n}$ are orthogonal matrices, and $\Sigma\in\mathbb{R}^{m\times n}$ is a diagonal matrix with real non-negative diagonal entries $\Sigma_{ii}=\sigma_i$, $i=1,\ldots,\min\{m,n\}$. The vectors $u_i$ and $v_i$ are called the left and right singular vectors, respectively, and the $\sigma_i$ are the singular values.

It is customary to write the decomposition in a way that the singular values are sorted in non-increasing order, $\sigma_1\geq\sigma_2\geq\ldots\geq\sigma_r>\sigma_{r+1}=\ldots=\sigma_n=0$, where $r=\mathrm{rank}(A)$. We can write the decomposition~\eqref{eq:svd} as a sum of outer product matrices,
\begin{equation}
\label{eq:svdouter}
A=\sum_{i=1}^{r}\sigma_iu_iv_i^T.
\end{equation}
It is well known that if only $k<r$ terms in~\eqref{eq:svdouter} are considered, the resulting matrix is the best rank-$k$ approximation of matrix $A$, in the least squares sense. This so called truncated SVD of $A$ is what one can usually afford to compute in the large-scale, sparse case. More generally, we will consider the case where the $k<r$ terms taken in~\eqref{eq:svdouter} correspond to either the largest or the smallest singular values, and we will refer to this decomposition as the partial SVD.

\begin{figure}
\centering
\begin{tikzpicture}[scale=0.5]
  \draw[fill=black!5] (-6.25,0) rectangle node {$U_A^T$} +(6,6);
  \draw[fill=black!5] (0,0) rectangle node {$A$} +(3,6);
  \draw[fill=black!5] (3.25,3) rectangle node {$G$} +(3,3);
  \node at (6.75,4.5) {=};
  \draw (7.25,0) rectangle +(3,6) (7.25,6) -- +(3,-3);
  \node[right] at (10.25,4.5) {$C$};
  \begin{scope}[yshift=-3.5cm]
  \draw[fill=black!5] (-2.25,1) rectangle node {$U_B^T$} +(2,2);
  \draw[fill=black!5] (0,1) rectangle node {$B$} +(3,2);
  \draw[fill=black!5] (3.25,0) rectangle node {$G$} +(3,3);
  \node at (6.75,2) {=};
  \draw (7.25,1) rectangle +(3,2) ++(3,0) -- +(-2,2);
  \node[right] at (10.25,2) {$S$};
  \end{scope}
\end{tikzpicture}
\caption{\label{fig:gsvd}Scheme of the GSVD of two matrices $A$ and $B$, for the case $m>n$ and $p<n$.}
\end{figure}
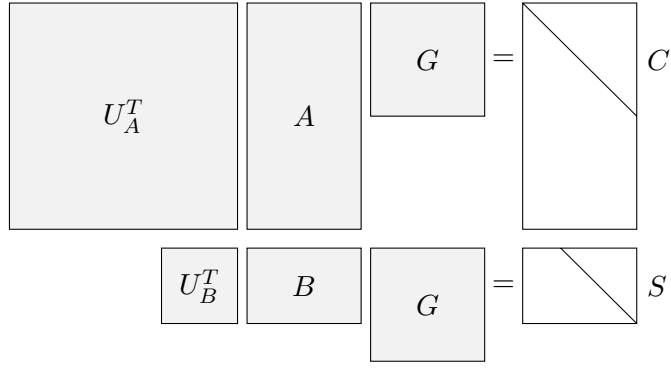

Now consider two matrices with the same column dimension, $A\in\mathbb{R}^{m\times n}$ and $B\in\mathbb{R}^{p\times n}$. The GSVD~\eqref{eq:gsvd} can also be written as
\begin{equation}
\label{eq:gsvd2}
A=U_ACG^{-1},\qquad B=U_BSG^{-1},
\end{equation}
with $U_A\in\mathbb{R}^{m\times m}$ and $U_B\in\mathbb{R}^{p\times p}$ orthogonal and $G\in\mathbb{R}^{n\times n}$ nonsingular. For our purpose, we can think of $C$ and $S$ as being diagonal matrices, but in the general case this needs a detailed discussion that we summarize below. In \eqref{eq:gsvd} and \eqref{eq:gsvd2}, we are assuming that the pair $\{A,B\}$ is regular, which means that the matrix obtained by stacking $A$ and $B$ has full column rank and hence the triangular factor of its QR decomposition is nonsingular\footnote{If the pair $\{A,B\}$ is not regular, a rank-revealing decomposition should be used instead of the QR~\cite{Paige:1981:TGS}.},
\begin{equation}
Z:=
\begin{bmatrix}
A\\
B
\end{bmatrix}
=
Q R
=
\begin{bmatrix}
Q_A\\
Q_B
\end{bmatrix}
R,
\label{eq-QR}
\end{equation}
where $R\in \mathbb{R}^{n\times n}$ is upper triangular and $Q\in\mathbb{R}^{(m+p)\times n}$ has orthonormal columns.
Then the structure of $C$ and $S$ is
\begin{equation}\label{eq:structure-cs}
\begin{array}{cc}
\begin{array}{c}
\begin{array}{ccc}
\!\!{\scriptstyle q} & \;{\scriptstyle \ell} & {\scriptstyle n-q-\ell}\!\!
\end{array}
\\
C=
\begin{bmatrix}
I_q & & \\
 & \hat{C} & \\
 & & O
\end{bmatrix}
\begin{array}{l}
{\scriptstyle q} \\ {\scriptstyle \ell} \\ {\scriptstyle m-q-\ell}
\end{array}
\end{array}
&
\begin{array}{c}
\begin{array}{ccc}
\!\!\!\!\!\!\!\!\!\!{\scriptstyle q} & \;\;{\scriptstyle \ell} & \;{\scriptstyle n-q-\ell}\!\!
\end{array}
\\
S=
\begin{bmatrix}
O & & \\
 & \hat{S} & \\
 & & I_{n-q-\ell}
\end{bmatrix}
\begin{array}{l}
{\scriptstyle p-n+q} \\ {\scriptstyle \ell} \\ {\scriptstyle n-q-\ell}
\end{array}
\end{array}
\end{array}
\end{equation}
where $\hat{C}$ and $\hat{S}$ are square diagonal matrices, $I_q$ and $I_{n-q-\ell}$ are identity matrices of the indicated size, and $O$ represents a rectangular block of zeros with the appropriate dimensions. Writing $\hat{C}=\mathrm{diag}(c_{q+1},\dots,c_{q+\ell})$ with $c_{q+1}\geq\cdots\geq c_{q+\ell}>0$, and $\hat{S}=\mathrm{diag}(s_{q+1},\dots,s_{q+\ell})$ with $0<s_{q+1}\leq\cdots\leq s_{q+\ell}$, we have that $c_i^2+s_i^2=1$, $i=q+1,\dots,q+\ell$, and the generalized singular values $\sigma(A,B)$ are the ratios of these quantities,
\begin{equation}\label{eq:xx}
\underbrace{\infty,\dots,\infty}_q,
\underbrace{c_{q+1}/s_{q+1},\dots,c_{q+\ell}/s_{q+\ell}}_\ell,
\underbrace{0,\dots,0}_{n-q-\ell}.
\end{equation}
In order to better understand the structure of $C$ and $S$, it is helpful to note that the zero blocks in~\eqref{eq:structure-cs} may have zero rows or columns, and also that the number of nonzero rows of $C$ is equal to $\mathrm{rank}(A)$ and similarly for $S$ with respect to $\mathrm{rank}(B)$. For instance, \cref{fig:gsvd} shows an example of GSVD for the case that $m>n$ and $p<n$. Assuming that $A$ and $B$ have full column and row ranks, respectively, we have that there are $q=n-p$ infinite generalized singular values.

As in the case of the SVD, for large-scale problems we will consider a partial GSVD, that is, we employ methods that compute approximations of $k$ quadruples $(\sigma_i,u^A_i,u^B_i,g_i)$ corresponding to either the largest or the smallest generalized singular values $\sigma_i$. We call $u^A_i$ and $u^B_i$ the left generalized singular vectors, while $g_i$ are the right generalized singular vectors. Note that if $\sigma_i$ are the generalized singular values of $\{A,B\}$, then $\sigma_i^{-1}$ are the generalized singular values of $\{B,A\}$.

\subsection{Equivalent eigenvalue problems}\label{sec:equivep}

The solution of the two problems presented in the previous section can be approached by formulating a related eigenvalue problem. More precisely, there are two possible strategies, that we will refer to as \emph{cross} and \emph{cyclic}. We start by discussing this in the context of the SVD and then extend it to the GSVD.

The SVD relation~\eqref{eq:svd} can be written as $AV=U\Sigma$ or as $A^TU=V\Sigma^T$. Suppose that $m\geq n$, then equating the columns we have
\begin{align}
\label{eq:svdleft}
Av_i&=u_i\sigma_i,\quad i=1,\ldots,n,\\
\label{eq:svdright}
A^Tu_i&=v_i\sigma_i,\quad i=1,\ldots,n,\\
\label{eq:svdright2}
A^Tu_i&=0,\quad i=n+1,\ldots,m.
\end{align}
Premultiplying \eqref{eq:svdleft} by $A^T$ and using \eqref{eq:svdright} results in the relation
\begin{equation}
\label{eq:eigleft}
A^TAv_i=\sigma_i^2v_i,
\end{equation}
that is, the $v_i$ are the eigenvectors of the symmetric matrix $A^TA$ corresponding to eigenvalues $\sigma_i^2$. If the corresponding left singular vectors are also required, they can be computed as $u_i=\frac{1}{\sigma_i}Av_i$ from \eqref{eq:svdleft}. Alternatively, it is possible to compute the left vectors first, via
\begin{equation}
\label{eq:eigright}
AA^Tu_i=\sigma_i^2u_i,
\end{equation}
and then the right ones as $v_i=\frac{1}{\sigma_i}A^Tu_i$, but care must be taken that~\eqref{eq:eigright} has at least $m-n$ zero eigenvalues. In practice, one would generally use~\eqref{eq:eigleft} if $m\geq n$ and~\eqref{eq:eigright} otherwise. We will call this approach the cross product eigenproblem.

The second strategy is the cyclic eigenproblem. Consider the symmetric matrix of order $m+n$
\begin{equation}
\label{eq:cyclic}
H(A)=\begin{bmatrix}0&A\\A^T&0\end{bmatrix},
\end{equation}
that has eigenvalues $\pm\sigma_i$, $i=1,\ldots,r$, together with $m+n-2r$ zero eigenvalues, where $r=\mathrm{rank}(A)$. The normalized eigenvectors corresponding to $\pm\sigma_i$ are $\frac{1}{\sqrt{2}}\left[\begin{smallmatrix}\pm u_i\\v_i\end{smallmatrix}\right]$. Hence we can extract the singular triplets $(\sigma_i,u_i,v_i)$ of $A$ directly from the eigenpairs of $H(A)$. Note that in this case the singular values are not squared, so the computed smallest singular values will not suffer from severe loss of accuracy as in the cross product approach. The drawback in this case is that small eigenvalues are located in the interior of the spectrum.

The cross and cyclic schemes can also be applied to the GSVD~\eqref{eq:gsvd}. The columns of $G$ satisfy
\begin{equation}
\label{eq:gsvdeigcross}
s_i^2A^TAg_i=c_i^2B^TBg_i,
\end{equation}
so solving a symmetric-definite generalized eigenvalue problem for the pencil $(A^TA,B^TB)$ provides us with the generalized singular values and right generalized singular vectors. From $g_i$ we can compute $u_i^A=\frac{1}{c_i}Ag_i$ and $u_i^B=\frac{1}{s_i}Bg_i$. This is the analog of the cross product eigenproblem for the SVD~\eqref{eq:eigleft}. It has the same concerns regarding the loss of accuracy in the smallest generalized singular values, but in this case one may consider computing these values as the reciprocals of the largest eigenvalues of the reversed pencil $(B^TB,A^TA)$.

Likewise, the formulation that is analogous to the eigenproblem associated with the cyclic matrix~\eqref{eq:cyclic} is to solve the symmetric-definite generalized eigenvalue problem defined by any of the matrix pencils
\begin{equation}
\label{eq:gsvdeigcyclic}
\left(
\begin{bmatrix}0&A\\A^T&0\end{bmatrix},
\begin{bmatrix}I&0\\0&B^TB\end{bmatrix}
\right),
\qquad\text{or}\qquad
\left(
\begin{bmatrix}0&B\\B^T&0\end{bmatrix},
\begin{bmatrix}I&0\\0&A^TA\end{bmatrix}
\right),
\end{equation}
of dimensions $m+n$ and $p+n$, respectively. The nonzero eigenvalues of the first pencil are $\pm\sigma_i$, while those of the second pencil are $\pm\sigma_i^{-1}$, so the latter is likely to be preferred when the smallest generalized singular values are required. The generalized singular vectors can be obtained from the eigenvectors corresponding to nonzero eigenvalues $\pm\sigma_i$ of the first pencil, $\frac{1}{\sqrt{2}}\left[\begin{smallmatrix}u_i^A\\\pm g_i/s_i\end{smallmatrix}\right]$, or the second pencil, $\frac{1}{\sqrt{2}}\left[\begin{smallmatrix}u_i^B\\\pm g_i/c_i\end{smallmatrix}\right]$. See~\cite{Huang:2020:CFC} for considerations on which of the two pencils is best in finite precision computations, based on the conditioning of $A$ and $B$.

Note that the symmetric-definite generalized eigenproblems discussed in this section may in fact be semi-definite, e.g., if $B^TB$ is singular because $p<n$. This may cause numerical difficulties when solving the problem via the cross and cyclic approaches.

\subsection{SVD via Lanczos bidiagonalization}\label{sec:lanczossvd}

In this section we give an overview of the computation of the partial SVD by means of Lanczos recurrences. Additional details can be found, e.g., in~\cite{Hernandez:2008:REP}. Without loss of generality, we will assume that $A$ is tall, i.e., $m\geq n$.

In the same way that the solution of the cross product eigenproblem for $A^TA$ can be approached by first computing an orthogonal similarity transformation to tridiagonal form, bidiagonalization methods for the SVD rely on an orthogonal reduction to bidiagonal form,
\begin{equation}\label{eq:fullbidiag}
A=P_nJ_nQ_n^T,
\end{equation}
with $P_n\in\mathbb{R}^{m\times n}$ and $Q_n\in\mathbb{R}^{n\times n}$ having orthonormal columns and $J_n\in\mathbb{R}^{n\times n}$ being upper bidiagonal, so that $J_n^TJ_n$ is tridiagonal with eigenvalues $\sigma_i^2$, that is, $J_n$ has the same singular values as $A$.

In a Lanczos method, we compute a partial bidiagonalization instead of the full one. From~\eqref{eq:fullbidiag} we can establish the two equalities $AQ_n=P_nJ_n$ and $A^TP_n=Q_nJ_n^T$, and equating the first $k<n$ columns we obtain the Lanczos relations
\begin{align}
\label{eq:lanczos-dec-1}
AQ_k &= P_kJ_k, \\
\label{eq:lanczos-dec-2}
A^TP_k &= Q_kJ_k^T+\beta_kq_{k+1}e_k^T,
\end{align}
where $J_k$ denotes the $k\times k$ leading principal submatrix of $J_n$, and we have used the notation
\begin{equation}
\label{eq:bidiagonal}
J_k=\begin{bmatrix}
\alpha_1 & \beta_1 \\
& \alpha_2 & \beta_2 \\
& & \ddots & \ddots \\
& & & \alpha_{k-1} & \beta_{k-1} \\
& & & & \alpha_k \\
\end{bmatrix}.
\end{equation}
Equating the $j$th column of~\eqref{eq:lanczos-dec-1}-\eqref{eq:lanczos-dec-2} gives the familiar double Lanczos recurrence that is shown in algorithmic form in \cref{alg:bidiagsvd}.

\newcommand{\normalize}[2]{Normalize: $#1=\|#2\|_2$, $#2=#2/#1$}

\begin{algorithm}
  \caption{Lanczos bidiagonalization}
  \label{alg:bidiagsvd}
  \begin{algorithmic}[1]
    \REQUIRE Matrix $A\in\mathbb{R}^{m\times n}$, unit-norm vector $q_1\in\mathbb{R}^n$, number of steps $k$.
    \ENSURE Partial bidiagonalization~\eqref{eq:lanczos-dec-1}-\eqref{eq:lanczos-dec-2}.
    \STATE Set $\beta_0=0$
    \FOR{$j=1,2,\dots,k$}
    \STATE \label{alg:bidiagsvd:expanda} $p_j=Aq_j-\beta_{j-1} p_{j-1}$
    \STATE \normalize{\alpha_j}{p_j}
    \STATE \label{alg:bidiagsvd:expandat} $q_{j+1}=A^Tp_j-\alpha_j q_j$
    \STATE \normalize{\beta_j}{q_{j+1}}
    \ENDFOR
  \end{algorithmic}
\end{algorithm}

It is possible to establish an equivalence between the output of \cref{alg:bidiagsvd} and the quantities computed by the Lanczos recurrence for tridiagonalizing the cross product matrix $A^TA$, see for instance~\cite{Hernandez:2008:REP}. In particular, the right Lanczos vectors $q_j$ computed by \cref{alg:bidiagsvd} form an orthonormal basis of the Krylov subspace $\mathcal{K}_k(A^TA,q_1)$. Similarly, the left Lanczos vectors $p_j$ span the Krylov subspace $\mathcal{K}_k(AA^T,Aq_1)$. Finally, there is also an equivalence with the Lanczos tridiagonalization associated with the cyclic matrix~\eqref{eq:cyclic}, provided that the initial vector $[0^T q_1^T]^T$ is used.

From the discussion above, it is clear that Ritz approximations of the singular triplets of $A$ can be obtained. After $k$ Lanczos steps, the Ritz values $\tilde{\sigma}_i$ (approximate singular values of $A$) are computed as the singular values of~$J_k$, and the Ritz vectors are $\tilde{u}_i=P_kx_i$ and $\tilde{v}_i=Q_ky_i$, where $x_i$ and $y_i$ are the left and right singular vectors of $J_k$, respectively.

\subsection{Residual and stopping criterion}\label{sec:residual}

As the number of Lanczos steps increase, the Ritz approximations become increasingly accurate. A criterion is required to determine when a certain approximate singular triplet $(\tilde{\sigma}_i,\tilde{u}_i,\tilde{v}_i)$ can be declared converged. And for this, we need to define a residual.

Due to the equivalence discussed in \cref{sec:equivep}, we can define the residual in terms of an equivalent eigenproblem. In the case of the cross product eigenproblem, \eqref{eq:eigleft} or \eqref{eq:eigright}, only one of the singular vectors would appear, so it is better to define the residual vector in terms of the eigenproblem associated with the cyclic matrix~\eqref{eq:cyclic}, whose norm is
\begin{equation}\label{eq:residualsvd}
\|r_i^\mathrm{SVD}\|_2=\sqrt{\|A\tilde{v}_i-\tilde{\sigma}_i\tilde{u}_i\|_2^2+\|A^T\tilde{u}_i-\tilde{\sigma}_i\tilde{v}_i\|_2^2}.
\end{equation}
This residual norm can be used in a posteriori error bounds. In particular, there exists a singular value $\sigma_{i'}$ of $A$ such that $|\sigma_{i'}-\tilde{\sigma}_i|\leq \frac{1}{\sqrt{2}}\|r_i^\mathrm{SVD}\|_2$~\cite{Bai:2000:TSA}.

In the context of Lanczos, the residual~\eqref{eq:residualsvd} can be computed cheaply as follows. If the Lanczos relations~\eqref{eq:lanczos-dec-1} and \eqref{eq:lanczos-dec-2} are multiplied on the right respectively by $y_i$ and $x_i$, the right and left singular vectors of $J_k$, then
$$
A\tilde{v}_i=\tilde{\sigma}_i\tilde{u}_i,\qquad
A^T\tilde{u}_i=\tilde{\sigma}_i\tilde{v}_i+\beta_kq_{k+1}e_k^Tx_i,
$$
and therefore
\begin{equation}
\label{eq:residual-beta}
\|r_i^\mathrm{SVD}\|_2=\beta_k|e_k^Tx_i|.
\end{equation}
The residual estimate~\eqref{eq:residual-beta} can be used in the stopping criterion in practical implementations of Lanczos bidiagonalization.

\subsection{Thick restart for the SVD}\label{sec:restart}

As in any Lanczos process, when run in finite precision arithmetic, loss of orthogonality among Lanczos vectors occurs in \cref{alg:bidiagsvd} whenever the Ritz values start to converge. The simplest cure is full reorthogonalization, which we consider in this work. This solution would be implemented by replacing \cref{alg:bidiagsvd:expanda} of \cref{alg:bidiagsvd} with
\begin{equation*}
p_j=\mathrm{orthog}(Aq_j,P_{j-1})
\end{equation*}
that applies Gram-Schmidt to explicitly orthogonalize vector $Aq_j$ against the columns of $P_{j-1}$, and similarly for \cref{alg:bidiagsvd:expandat}.

Full reorthogonalization requires keeping all previously computed (left and right) Lanczos vectors, and this justifies the need of a restart technique that limits the number of vectors, not only due to storage requirements but also because the computational cost of full reorthogonalization is proportional to the number of involved vectors.

Thick restart is very effective compared to explicit restart, because instead of explicitly building a new initial vector to rerun the algorithm, it keeps a smaller dimensional subspace that retains the most relevant spectral information computed so far. The computation is thus a sequence of subspace expansions and contractions, until the subspace contains enough converged solutions. The key is to purge unwanted components during the contraction, which is beneficial for overall convergence.

Suppose we have built the Lanczos relations~\eqref{eq:lanczos-dec-1}-\eqref{eq:lanczos-dec-2} of order $k$ and we want to compress to size $r<k$. We start by transforming the decomposition in a way that approximate singular values and the residual norms~\eqref{eq:residual-beta} appear explicitly in the equations, as follows. First, compute the SVD of the bidiagonal matrix, $J_k=X_k\tilde{\Sigma}_kY_k^T$, with $\tilde{\Sigma}_k=\mathrm{diag}(\tilde{\sigma}_1,\dots,\tilde{\sigma}_k)$. Post-multiplication of~\eqref{eq:lanczos-dec-1}-\eqref{eq:lanczos-dec-2} by $Y_k$ and $X_k$, respectively, results in
\begin{align}
\label{eq:svd-compress-1}
A\tilde{V}_k &= \tilde{U}_k\tilde{\Sigma}_k, \\
\label{eq:svd-compress-2}
A^T\tilde{U}_k &= \tilde{V}_k\tilde{\Sigma}_k+\beta_kq_{k+1}e_k^TX_k,
\end{align}
where $\tilde{U}_k=P_kX_k$ and $\tilde{V}_k=Q_kY_k$, whose columns are the left and right Ritz vectors. This would be an exact partial SVD of $A$ if it were not for the second summand of the right-hand side of~\eqref{eq:svd-compress-2}, which is related to the residual norms of the $k$ Ritz approximations. Equation~\eqref{eq:svd-compress-2} can also be written as
\begin{equation}\label{eq:svd-compress-2-alt}
A^T\tilde{U}_k = \begin{bmatrix}\tilde{V}_k & q_{k+1}\end{bmatrix}\tilde{J}_k^T,\qquad
\tilde{J}_k=\begin{bmatrix}
  \tilde{\sigma}_1 & & & & \rho_1\\
  & \tilde{\sigma}_2 & & & \rho_2\\
  & & \ddots & & \vdots\\
  & & & \tilde{\sigma}_k & \rho_k\\
\end{bmatrix},
\end{equation}
where $\rho_i:=\beta_ke_k^Tx_i$. Note that $|\rho_i|$ is equal to the residual norm~\eqref{eq:residual-beta}. Equation~\eqref{eq:svd-compress-2-alt} is represented graphically in \cref{fig:thick-svd} (step 2).

\def\mydiagonal{ 
  foreach \ii in {1,...,\nits} { -- ++(0,-0.2) -- ++(0.2,0) }
  foreach \ii in {1,...,\nits} { -- ++(0,0.2) -- ++(-0.2,0) }
}
\def\myupperbidiagonal{ 
  foreach \ii in {1,...,\nits} { -- ++(0,-0.2) -- ++(0.2,0) }
  -- ++(0,0.2)
  foreach \ii in {1,...,\nitsmone} { -- ++(0,0.2) -- ++(-0.2,0) }
  -- ++(-0.2,0)
}
\def\myupperbidiagonalextra{ 
  foreach \ii in {1,...,\nits} { -- ++(0,-0.2) -- ++(0.2,0) }
  -- ++(0.2,0)
  foreach \ii in {1,...,\nits} { -- ++(0,0.2) -- ++(-0.2,0) }
  -- ++(-0.2,0)
}
\def\mylowerbidiagonal{ 
  -- ++(0,-0.2)
  foreach \ii in {1,...,\nitsmone} { -- ++(0,-0.2) -- ++(0.2,0) }
  -- ++(0.2,0)
  foreach \ii in {1,...,\nits} { -- ++(0,0.2) -- ++(-0.2,0) }
}
\def\mylowerbidiagonalextra{ 
  -- ++(0,-0.2)
  foreach \ii in {1,...,\nits} { -- ++(0,-0.2) -- ++(0.2,0) }
  -- ++(0,0.2)
  foreach \ii in {1,...,\nits} { -- ++(0,0.2) -- ++(-0.2,0) }
}

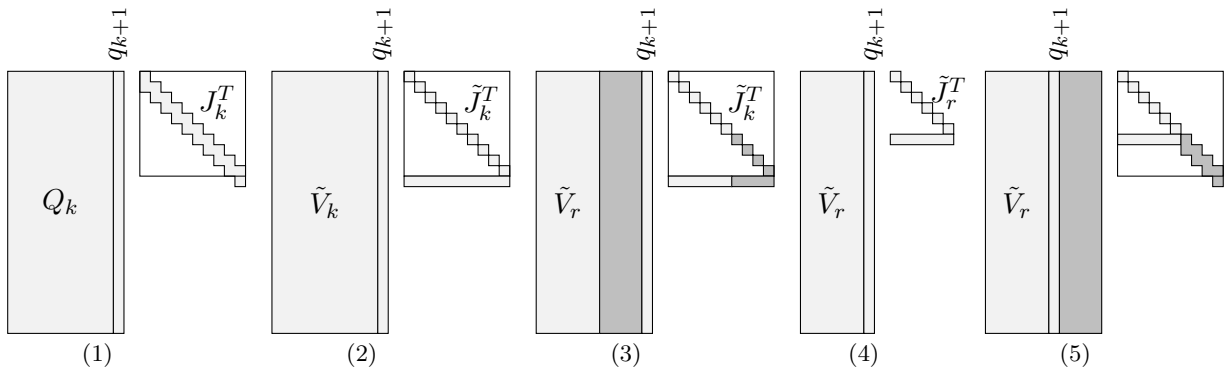
\begin{figure}
  \resizebox{\textwidth}{!}{\begin{tikzpicture}[scale=0.7]
    \begin{scope}[xshift=0]
      \draw[fill=black!5] (8.3,0) rectangle node {$Q_k$} +(2,5);
      \draw[fill=black!5] (10.3,0) rectangle +(0.2,5);
      \node[right,rotate=90] at (10.4,5) {$q_{k+1}$};
      \def\nits{10}
      \draw[fill=black!5] (10.8,5) \mylowerbidiagonalextra;
      \draw (10.8,5) rectangle ++(0.2*\nits,-0.2*\nits);
      \node[left] at (12.8,4.4) {$J_k^T$};
      \node[below,font=\footnotesize] at (10,0) {(1)};
    \end{scope}
    \begin{scope}[xshift=5cm]
      \draw[fill=black!5] (10.3,0) rectangle +(0.2,5);
      \node[right,rotate=90] at (10.4,5) {$q_{k+1}$};
      \draw[fill=black!5] (8.3,0) rectangle node {$\tilde{V}_k$} +(2,5);
      \def\nits{10}
      \draw[fill=black!5] (10.8,5) \mydiagonal;
      \draw (10.8,5) rectangle ++(0.2*\nits,-0.2*\nits);
      \node[left] at (12.8,4.4) {$\tilde{J}_k^T$};
      \draw[fill=black!5] (10.8,2.8) rectangle +(2,0.2);
      \node[below,font=\footnotesize] at (10,0) {(2)};
    \end{scope}
    \begin{scope}[xshift=10cm]
      \draw[fill=black!5] (8.3,0) rectangle node {$\tilde{V}_r$} +(1.2,5);
      \draw[fill=black!25] (8.3,0) ++(1.2,0) rectangle +(0.8,5);
      \draw[fill=black!5] (10.3,0) rectangle +(0.2,5);
      \node[right,rotate=90] at (10.4,5) {$q_{k+1}$};
      \def\nits{6}
      \draw[fill=black!5] (10.8,5) \mydiagonal;
      \def\nits{4}
      \draw[fill=black!25] (12,3.8) \mydiagonal;
      \draw (10.8,5) rectangle ++(0.2*10,-0.2*10);
      \node[left] at (12.8,4.4) {$\tilde{J}_k^T$};
      \draw[fill=black!5] (10.8,2.8) rectangle +(1.2,0.2);
      \draw[fill=black!25] (12,2.8) rectangle +(0.8,0.2);
      \node[below,font=\footnotesize] at (10,0) {(3)};
    \end{scope}
    \begin{scope}[xshift=15cm]
      \draw[fill=black!5] (8.3,0) rectangle node {$\tilde{V}_r$} +(1.2,5);
      \draw[fill=black!5] (9.5,0) rectangle +(0.2,5);
      \node[right,rotate=90] at (9.7,5) {$q_{k+1}$};
      \def\nits{6}
      \draw[fill=black!5] (10,5) \mydiagonal;
      \node at (11.1,4.6) {$\tilde{J}_r^T$};
      \draw[fill=black!5] (10,3.6) rectangle +(1.2,0.2);
      \node[below,font=\footnotesize] at (9.5,0) {(4)};
    \end{scope}
    \begin{scope}[xshift=18.5cm]
      \draw[fill=black!5] (8.3,0) rectangle node {$\tilde{V}_r$} +(1.2,5);
      \draw[fill=black!5] (9.5,0) rectangle +(0.2,5);
      \node[right,rotate=90] at (9.7,5) {$q_{k+1}$};
      \def\nits{6}
      \draw[fill=black!5] (10.8,5) \mydiagonal;
      \draw[fill=black!5] (10.8,3.6) rectangle +(1.2,0.2);
      \draw[fill=black!25] (8.3,0) ++(1.4,0) rectangle +(0.8,5);
      \def\nits{4}
      \draw[fill=black!25] (12,3.8) \mylowerbidiagonalextra;
      \draw (10.8,5) rectangle ++(0.2*10,-0.2*10);
      \node[below,font=\footnotesize] at (10,0) {(5)};
    \end{scope}
  \end{tikzpicture}}
  \caption{Illustration of the five steps of thick restart in SVD: (1) initial Lanczos factorization of order $k$, (2) solve projected problem, (3) sort and check convergence, (4) truncate to factorization of order $r$, and (5) extend to a factorization of order $k$.}
  \label{fig:thick-svd}
\end{figure}

Due to the fact that $\tilde{\Sigma}_k$ in~\eqref{eq:svd-compress-1}-\eqref{eq:svd-compress-2} is diagonal, it is possible to truncate this decomposition at any size $r$. For this, we must permute it so that the leading principal submatrix of $\tilde{\Sigma}_k$ contains the approximations of the wanted singular values (typically the largest or smallest ones). \Cref{fig:thick-svd} (step 3) shows in dark gray the part of the decomposition that will be discarded. Then the new decomposition after truncation is
\begin{align*}
A\tilde{V}_r &= \tilde{U}_r\tilde{\Sigma}_r, \\
A^T\tilde{U}_r &= \tilde{V}_r\tilde{\Sigma}_r+q_{k+1}b_r^T,
\end{align*}
where $b_r^T=[\rho_1,\dots,\rho_r]$, see \cref{fig:thick-svd} (step 4). It only remains to extend the decomposition by running a modified version of \cref{alg:bidiagsvd} that starts the loop at $j=r+1$. Note that the first newly generated left Lanczos vector is computed as $p_{r+1}=\mathrm{orthog}(Aq_{k+1},\tilde{U}_r)$, whose orthogonalization coefficients are precisely the $\rho_i$'s. When the algorithm stops at iteration $j=k$, the new $J_k$ matrix has the shape depicted in \cref{fig:thick-svd} (step 5), a bidiagonal except for the leading part that has an arrowhead form.

\subsection{GSVD via joint Lanczos bidiagonalization}\label{sec:lanczosgsvd}

We now turn our attention to the GSVD, and consider Lanczos methods that compute a joint bidiagonalization of the two matrices, $A$ and $B$. In this context, the stacked matrix $Z$~\eqref{eq-QR} is relevant and will appear in the algorithms. Also, the matrices $Q_A$ and $Q_B$ in~\eqref{eq-QR}, whose row dimensions are the same as $A$ and $B$, respectively, will be used for deriving the methods, but need not be formed explicitly as justified later.

The matrices $C$ and $S$~\eqref{eq:structure-cs} are related to the CS decomposition~\cite[\S 2.6.4]{Golub:1996:MC} of $\{Q_A,Q_B\}$, that is, the values $c_i$ and $s_i$ are related to the singular values of $Q_A$ and $Q_B$, respectively. Let the CS decomposition of $\{Q_A,Q_B\}$ be
\begin{equation}
Q_A = U_A C W^T,\quad Q_B = U_B S W^T,
\label{eq-Q-CS-decomp}
\end{equation}
where $U_A$, $U_B$ are the same as in~\eqref{eq:gsvd2}, and $W\in\mathbb{R}^{n\times n}$ is also orthogonal. Note that the first equation in~\eqref{eq-Q-CS-decomp} can be seen as the conventional SVD of $Q_A$, with the singular values sorted in non-increasing order, while the second equation is an SVD-like relation where the singular values appear in non-decreasing order with the largest value at the bottom-right corner of $S$.

If $Z$ has full column rank, the GSVD decomposition of $\{A,B\}$ is given by
\begin{equation}
A = U_A C G^{-1},\quad B = U_B S G^{-1},
\label{eq-gsvd}
\end{equation}
where $G = R^{-1} W$ with $R$ as defined in~\eqref{eq-QR}.

An approach to compute the CS decomposition of $\{Q_A,Q_B\}$ is to perform a bidiagonalization of both matrices, i.e., a decomposition of the form
\begin{equation}
Q_A = U J V^T,\quad Q_B = \widehat{U} \widehat{J} V^T,
\label{eq-Q-bidiag}
\end{equation}
with $J$, $\widehat{J}$ (upper or lower) bidiagonal matrices such that the stacked matrix~$\left[\begin{smallmatrix}J\\\widehat{J}\end{smallmatrix}\right]$ has orthonormal columns, that is, $J^T J + \widehat{J}^T \widehat{J} = I$. Then, the problem is reduced to the CS decomposition of the bidiagonal matrices $\{J, \widehat{J}\}$,
\begin{equation}
J = X C Y^T,\quad \widehat{J} = \widehat{X} S Y^T.
\label{eq-J-CS-decomp-full-size}
\end{equation}
Substituting in~\eqref{eq-Q-bidiag} we get the CS decomposition~\eqref{eq-Q-CS-decomp} with $U_A = U X$, $U_B = \widehat{U} \widehat{X}$ and $W = V Y$.

As mentioned in \cref{sec:svd}, if $A$, $B$ are very large and sparse matrices, the interest is normally to compute a partial decomposition, i.e., a few (extreme) generalized singular values and vectors. In that case, a partial bidiagonalization of matrices $Q_A$ and $Q_B$ is done using Lanczos recurrences, as described next.

Zha~\cite{Zha:1996:CGS} presented an algorithm for the joint bidiagonalization of $Q_A$ and $Q_B$, in which both matrices are reduced to upper bidiagonal form without explicitly computing $Q_A$ or $Q_B$.
Later, Kilmer \emph{et al.}~\cite{Kilmer:2007:PAG} proposed a variation of the joint bidiagonalization where $Q_A$ and $Q_B$ are transformed to lower and upper bidiagonal forms, respectively. To keep the presentation short, we focus on the latter variant from now on, although our solver also includes an implementation of Zha's variant.

Although not computing $Q_A$ or $Q_B$ explicitly, Kilmer's joint bidiagonalization is based on the application of the lower and upper Lanczos bidiagonalization algorithms in~\cite{Paige:1982:LAS} to $Q_A$ and $Q_B$, respectively, yielding
\begin{align}
Q_A V_k &= U_{k+1} J_k, &
Q_A^T U_{k+1} &= V_k J_k^T + \alpha_{k+1} v_{k+1} e_{k+1}^T,
\label{eq-kilmer-A}
\\
Q_B \widehat{V}_k &= \widehat{U}_{k} \widehat{J}_k, &
Q_B^T \widehat{U}_{k} &= \widehat{V}_k \widehat{J}_k^T + \hat{\beta}_{k} \hat{v}_{k+1} e_{k}^T,
\label{eq-kilmer-B}
\end{align}
with the column-orthonormal matrices $U_{k+1}=\left[u_1,u_2,\dots,u_{k+1}\right]$, $V_k=\left[v_1,v_2,\dots,v_{k}\right]$, $\widehat{U}_{k}=\left[\hat{u}_1,\hat{u}_2,\dots,\hat{u}_{k}\right]$ and $\widehat{V}_k=\left[\hat{v}_1,\hat{v}_2,\dots,\hat{v}_{k}\right]$ and the lower and upper bidiagonal matrices
$$
J_k=\begin{bmatrix}
\alpha_1 &         &           &\\
\beta_2  &\alpha_2 &           &\\
         & \ddots  & \ddots    &\\
         &         & \beta_{k} & \alpha_k\\
         &         &           & \beta_{k+1}
\end{bmatrix}
\in \mathbb{R}^{(k+1)\times k},
\quad
\widehat{J}_k=\begin{bmatrix}
\hat{\alpha}_1 & \hat{\beta}_1 &         &\\
               &\hat{\alpha}_2 & \ddots  & \\
               &               & \ddots  & \hat{\beta}_{k-1} \\
               &               &         & \hat{\alpha}_k\\
\end{bmatrix}
\in \mathbb{R}^{k\times k}.
$$
Note that the bottom Lanczos relations~\eqref{eq-kilmer-B} are analogous to those used for the bidiagonalization in SVD~\eqref{eq:lanczos-dec-1}-\eqref{eq:lanczos-dec-2}. However, the top Lanczos relations~\eqref{eq-kilmer-A} differ in that the associated bidiagonal matrix $J_k$ is lower bidiagonal with one more row than columns, and the basis of left Lanczos vectors $U_{k+1}$ contains one more vector than in the other case. The reason is that the lower bidiagonalization algorithm in~\cite{Paige:1982:LAS} starts the recurrence with the left vectors instead of the right ones.

Zha's method generates two upper bidiagonal matrices by applying \cref{alg:bidiagsvd} to both $Q_A$ and $Q_B$ with the same initial vector. In contrast, Kilmer's method uses the two types of bidiagonalizations and connects them by using the first right Lanczos vector $v_1$ generated in~\eqref{eq-kilmer-A} as initial vector $\hat{v}_1$ in~\eqref{eq-kilmer-B}.

It can be shown~\cite{Zha:1996:CGS,Kilmer:2007:PAG} that if $v_1=\hat{v}_1$, the joint bidiagonalization given by \eqref{eq-kilmer-A}-\eqref{eq-kilmer-B} verifies
\begin{equation}
\hat{v}_{i}=(-1)^{i-1} v_{i},\quad \hat{\alpha}_i\hat{\beta}_i=\alpha_{i+1}\beta_{i+1},\quad i=1,2,\dots.
\label{eq-relation-v}
\end{equation}
Thus, \eqref{eq-kilmer-A}-\eqref{eq-kilmer-B} can be rewritten as
\begin{align}
Q_A V_k &= U_{k+1} J_k, &
Q_A^T U_{k+1} &= V_k J_k^T + \alpha_{k+1} v_{k+1} e_{k+1}^T,
\label{eq-kilmer2-A}
\\
Q_B V_k &= \widehat{U}_{k} \check{J}_k, &
Q_B^T \widehat{U}_{k} &= V_k \check{J}_k^T + \check{\beta}_{k} v_{k+1} e_{k}^T,
\label{eq-kilmer2-B}
\end{align}
where $\check{J}_k=\widehat{J}_k D$, $D=\mathrm{diag}(1,-1,\dots,(-1)^{k-1})$ and $\check{\beta}_k=(-1)^k\hat{\beta}_{k}$.

Taking into account that $Q_A=\begin{bmatrix}I_m&0\end{bmatrix}Q$, $Q_B=\begin{bmatrix}0&I_p\end{bmatrix}Q$, and premultiplying equalities on the right side of \eqref{eq-kilmer2-A}-\eqref{eq-kilmer2-B} by $Q$, we get
\begin{align*}
\begin{bmatrix}I_m & 0\end{bmatrix} Q V_k &= U_{k+1} J_k, &
Q Q^T \begin{bmatrix}U_{k+1}\\0\end{bmatrix} &= Q V_k J_k^T + \alpha_{k+1} Q v_{k+1} e_{k+1}^T,
\\
\begin{bmatrix}0 & I_p\end{bmatrix}Q V_k &= \widehat{U}_{k} \check{J}_k, &
Q Q^T \begin{bmatrix}0\\\widehat{U}_{k}\end{bmatrix} &= Q V_k \check{J}_k^T + \check{\beta}_{k} Q v_{k+1} e_{k}^T,
\end{align*}
or, defining $\widetilde{V}_k=Q V_k$,
\begin{align}
\begin{bmatrix}I_m & 0\end{bmatrix} \widetilde{V}_k &= U_{k+1} J_k, &
Q Q^T \begin{bmatrix}U_{k+1}\\0\end{bmatrix} &= \widetilde{V}_k J_k^T + \alpha_{k+1} \tilde{v}_{k+1} e_{k+1}^T,
\label{eq-kilmer3-A}
\\
\begin{bmatrix}0 & I_p\end{bmatrix}\widetilde{V}_k &= \widehat{U}_{k} \check{J}_k, &
Q Q^T \begin{bmatrix}0\\\widehat{U}_{k}\end{bmatrix} &= \widetilde{V}_k \check{J}_k^T + \check{\beta}_{k} \tilde{v}_{k+1} e_{k}^T.
\label{eq-kilmer3-B}
\end{align}

The joint bidiagonalization process in~\cite{Kilmer:2007:PAG} uses the two Lanczos relations in~\eqref{eq-kilmer3-A} and the first one in~\eqref{eq-kilmer3-B} to compute matrices $U_{k+1}$, $\widehat{U}_k$, $\widetilde{V}_k$, $J_k$, $\check{J}_k$. Recall that the vectors $u_{j}$, $\hat{u}_j$, $\tilde{v}_j$ have lengths $m$, $p$ and $m+p$, respectively.
\Cref{alg:kilmer} gives the details of Kilmer's joint bidiagonalization, where in \cref{alg:kilmer:expand1,alg:kilmer:expand2} $\mathrm{expand}(A,B,u_{j+1})$ denotes the operation that generates new Krylov directions from the $A$ and $B$ matrices as follows.
Each step $j$ of the joint bidiagonalization requires computing $Q Q^T \tilde{u}_{j+1}$, where $\tilde{u}_{j+1}=\left[u_{j+1}^T,0\right]^T$. Note that this is the orthogonal projection of $\tilde{u}_{j+1}$ onto the column space of $Z$, which means that
$Q Q^T \tilde{u}_{j+1}=Z x_{j+1}$, where $x_{j+1}$ is the solution of the least squares problem
\begin{equation}\label{eq:leastsquares}
x_{j+1}=\underset{x\in\mathbb{R}^n}{\arg \min} \| Z x-\tilde{u}_{j+1}\|.
\end{equation}
The $\mathrm{expand}(A,B,u_{j+1})$ operation first computes the solution of the least squares problem~\eqref{eq:leastsquares} by padding $u_{j+1}$ below with $p$ zeros, and then performs an additional multiplication by $Z$.
If $Z$ is a large sparse matrix, the least squares problem is solved by means of an iterative solver such as the LSQR algorithm~\cite{Paige:1982:LAS}. Thus, the bidiagonalization process does not need to compute the QR factorization of $Z$.

\refstepcounter{myalg}
\begin{algorithm}
  \caption{Lower-upper joint Lanczos bidiagonalization~\cite{Kilmer:2007:PAG}}
  \label{alg:kilmer}
  \begin{algorithmic}[1]
    \REQUIRE Matrices $A\in\mathbb{R}^{m\times n}$, $B\in\mathbb{R}^{p\times n}$, unit-norm vector $u_1\in\mathbb{R}^m$, number of steps $k$.
    \ENSURE Partial joint bidiagonalization~\eqref{eq-kilmer3-A}-\eqref{eq-kilmer3-B}.
    \STATE Set $\hat{\beta}_0=0$
    \STATE \label{alg:kilmer:expand1} $\tilde{v}_1=\mathrm{expand}(A,B,u_1)$
    \STATE \normalize{\alpha_1}{\tilde{v}_1}
    \FOR{$j=1,2,\dots,k$}
    \STATE \label{alg:kilmer:matlab1} $\hat{u}_j=(-1)^{j-1}\left[0,I_p\right]\tilde{v}_j-\hat{\beta}_{j-1}\hat{u}_{j-1}$
    \STATE \normalize{\hat{\alpha}_j}{\hat{u}_j}
    \STATE \label{alg:kilmer:matlab2} $u_{j+1}=\left[I_m,0\right]\tilde{v}_j-\alpha_j u_j$
    \STATE \normalize{\beta_{j+1}}{u_{j+1}}
    \STATE \label{alg:kilmer:expand2} $\tilde{v}_{j+1}=\mathrm{expand}(A,B,u_{j+1})-\beta_{j+1}\tilde{v}_j$
    \STATE \normalize{\alpha_{j+1}}{\tilde{v}_{j+1}}
    \STATE \label{alg:kilmer:last-in-for} $\hat{\beta}_j=(\alpha_{j+1}\beta_{j+1})/\hat{\alpha}_j$
    \ENDFOR
  \end{algorithmic}
\end{algorithm}

After running \cref{alg:kilmer} we get matrices $U_{k+1}, \widehat{U}_k,\widetilde{V}_k,J_k,\check{J}_k$ that, taking $V_k=Q^T \widetilde{V}_k$, satisfy equations~\eqref{eq-kilmer2-A}-\eqref{eq-kilmer2-B}.
On the other hand, as pointed out in~\cite{Kilmer:2007:PAG}, it can be shown that
$J_k^T J_k+\check{J}_k^T \check{J}_k=I_k$, which means that the matrix pair $\lbrace J_k,\check{J}_k \rbrace$ admits a CS decomposition
\begin{equation}
J_k=X_{k+1}\begin{bmatrix}C_k\\0\end{bmatrix}Y_k^T,\qquad \check{J}_k=\widehat{X}_k S_k Y_k^T,
\label{eq-J-CS-decomp}
\end{equation}
where $X_{k+1} \in \mathbb{R}^{(k+1)\times (k+1)}$, $\widehat{X}_k,Y_k\in \mathbb{R}^{k\times k}$ are orthogonal matrices, and $C_k,S_k$ are $k\times k$ diagonal matrices. Note that these matrices do not correspond to subblocks of the matrices $X,\widehat{X},Y,C,S$ appearing in \eqref{eq-J-CS-decomp-full-size}.

Using the CS decomposition \eqref{eq-J-CS-decomp}, the Lanczos relations~\eqref{eq-kilmer2-A}-\eqref{eq-kilmer2-B} become
\begin{align}
Q_A V_k Y_k &= U_{k+1} X_{k+1} \begin{bmatrix}C_k\\0\end{bmatrix}, &
Q_A^T U_{k+1} X_{k+1} &= V_k Y_k \begin{bmatrix}C_k & 0\end{bmatrix} + \alpha_{k+1} v_{k+1} e_{k+1}^T X_{k+1},
\label{eq-kilmer-after-CSA}
\\
Q_B V_k Y_k &= \widehat{U}_{k} \widehat{X}_k S_k, &
Q_B^T \widehat{U}_{k} \widehat{X}_k &= V_k Y_k S_k + \check{\beta}_{k} v_{k+1} e_{k}^T \widehat{X}_k.
\label{eq-kilmer-after-CSB}
\end{align}
Taking $\tilde{w}_i$, $\tilde{u}_i^A$, $\tilde{u}_i^B$, $x_i$, $\hat{x}_i$ and $y_i$ as the $i$th columns of $V_k Y_k$, $U_{k+1} X_{k+1}$, $\widehat{U}_{k} \widehat{X}_k$, $X_k$, $\widehat{X}_k$ and $Y_k$, respectively, and $\tilde{c}_i$, $\tilde{s}_i$, as the $i$th diagonal elements of $C_k$, $S_k$, we have
\begin{align}
Q_A \tilde{w}_i &= \tilde{c}_i \tilde{u}_i^A, &
Q_A^T \tilde{u}_i^A &= \tilde{c}_i \tilde{w}_i + \alpha_{k+1} v_{k+1} e_{k+1}^T x_i,
\label{eq-kilmer-vectors-QA}
\\
Q_B \tilde{w}_i &= \tilde{s}_i \tilde{u}_i^B, &
Q_B^T \tilde{u}_i^B &= \tilde{s}_i \tilde{w}_i + \check{\beta}_{k} v_{k+1} e_{k}^T \hat{x}_i,
\label{eq-kilmer-vectors-QB}
\end{align}
for $i=1,2,\dots,k$.
Thus, $\tilde{u}_i^A$, $\tilde{u}_i^B$, are approximations of the left generalized singular vectors for $A$ and $B$, respectively, while the approximate right singular vectors $\tilde{g}_i=R^{-1}\tilde{w}_i$ can be computed as the solution to the least squares problem $Z \tilde{g}_i=\hat{w}_i$, where $\hat{w}_i=Q\tilde{w}_i=Q V_k y_i=\widetilde{V}_k y_i$.

\section{Thick restart for the GSVD}\label{sec:trlanczos}

With all the ingredients presented in the previous section, we are now able to adapt the thick restart technique of \cref{sec:restart} to the case of Kilmer's joint bidiagonalization for the GSVD. As was done for the SVD, the goal is to successively compress and expand the decomposition until the working Lanczos bases contain sufficiently good approximations to the wanted solutions. The compression is carried out by transforming the decomposition computed by \cref{alg:kilmer} to a form where the generalized singular values (or, more precisely, the $c_i$ and $s_i$ values) appear explicitly, and then truncating it in a way that keeps the most relevant part. This requires computing a small-sized GSVD (or CS decomposition), as in~\eqref{eq-J-CS-decomp}. Once the decomposition is truncated, it should be possible to extend it again by a slightly modified variant of \cref{alg:kilmer} whose loop starts at the current size after truncation.

\subsection{Restarting Kilmer's joint bidiagonalization}\label{sec:trkilmer}

Substituting \eqref{eq-J-CS-decomp} in \eqref{eq-kilmer3-A}-\eqref{eq-kilmer3-B}, we have
\begin{align*}
\begin{bmatrix}I_m& 0\end{bmatrix} \widetilde{V}_k Y_k &= U_{k+1} X_{k+1} \begin{bmatrix}C_k\\0\end{bmatrix},
&
Q Q^T \begin{bmatrix}U_{k+1} X_{k+1}\\0\end{bmatrix}
&= \widetilde{V}_k Y_k \begin{bmatrix}C_k&0\end{bmatrix} +
\alpha_{k+1} \tilde{v}_{k+1} e_{k+1}^T X_{k+1},
\\
\begin{bmatrix}0 & I_p\end{bmatrix} \widetilde{V}_k Y_k &= \widehat{U}_{k} \widehat{X}_k S_k,
&
Q Q^T \begin{bmatrix}0\\\widehat{U}_{k} \widehat{X}_k\end{bmatrix}
&= \widetilde{V}_k Y_k S_k + \check{\beta}_{k} \tilde{v}_{k+1} e_{k}^T \widehat{X}_k.
\end{align*}
This decomposition is similar to~\eqref{eq-kilmer-after-CSA}-\eqref{eq-kilmer-after-CSB}, but expressed in terms of the $\widetilde{V}_k$ basis instead of $V_k$.

In order to truncate the equations above, we define $Y_r$ and $\widehat{X}_r$ as the matrices formed by taking the first $r$ columns of $Y_k$ and $\widehat{X}_k$, respectively, $C_r$ and $S_r$ as the leading principal $r\times r$ submatrix of $C_k$ and $S_k$, and $X_{r+1}=[x_1,x_2\dots,x_r, x_{k+1}]$, where $x_j$ is column $j$ of $X_{k+1}$. Note that we keep the last column $x_{k+1}$ and append it as column $r+1$. We can then write
\begin{align*}
\begin{bmatrix}I_m& 0\end{bmatrix} \widetilde{V}_k Y_r &= U_{k+1} X_{r+1} \begin{bmatrix}C_r\\0\end{bmatrix},
&
Q Q^T \begin{bmatrix}U_{k+1}X_{r+1}\\0\end{bmatrix}
&= \widetilde{V}_k Y_r \begin{bmatrix}C_r&0\end{bmatrix} + \alpha_{k+1} \tilde{v}_{k+1} e_{k+1}^T X_{r+1},
\\
\begin{bmatrix}0 & I_p\end{bmatrix} \widetilde{V}_k Y_r &= \widehat{U}_{k} \widehat{X}_r S_r,
&
Q Q^T \begin{bmatrix}0\\\widehat{U}_{k}\widehat{X}_r\end{bmatrix}
&= \widetilde{V}_k Y_r S_r + \check{\beta}_{k} \tilde{v}_{k+1} e_{k}^T \widehat{X}_r,
\end{align*}
or
\begin{align}
\begin{bmatrix}I_m& 0\end{bmatrix} \widetilde{V}'_r &= U'_{r+1} \begin{bmatrix}C_r\\0\end{bmatrix},
&
Q Q^T \begin{bmatrix}U'_{r+1}\\0\end{bmatrix} &= \widetilde{V}'_r
\begin{bmatrix}C_r&0\end{bmatrix} +\tilde{v}_{k+1} b_{r+1}^T,
\label{eq-restart-a}
\\
\begin{bmatrix}0 & I_p\end{bmatrix} \widetilde{V}'_r &= \widehat{U}'_{r} S_r,
&
Q Q^T \begin{bmatrix}0\\\widehat{U}'_{r}\end{bmatrix}
&= \widetilde{V}'_r S_r + \tilde{v}_{k+1}\hat{b}_r^T,
\label{eq-restart-b}
\end{align}
where $b_{r+1}=\alpha_{k+1} X_{r+1}^T e_{k+1}$ and $\hat{b}_r=\check{\beta}_{k} \widehat{X}_r^T e_{k}$, while the orthonormal vector bases have been updated according to $U'_{r+1}=U_{k+1} X_{r+1}$, $\widehat{U}'_{r}=\widehat{U}_{k} \widehat{X}_{r}$ and $\widetilde{V}'_{r}=\widetilde{V}_{k} Y_{r}$.

The bidiagonalization process continues with the updated vector bases and taking $\tilde{v}_{k+1}$ as the next vector of the $\widetilde{V}$ basis. That is, \cref{alg:kilmer} is run with the loop starting at $j=r+1$.

\begin{figure}
  \resizebox{\textwidth}{!}{\begin{tikzpicture}[scale=0.85]
    \begin{scope}[xshift=0]
      \draw[fill=black!5] (8.3,0) rectangle node {$\widetilde{V}_k$} +(2,5);
      \draw[fill=black!5] (10.3,0) rectangle +(0.2,5);
      \node[right,rotate=90] at (10.4,5) {$\tilde{v}_{k+1}$};
      \def\nits{11}
      \def\nitsmone{10}
      \draw[fill=black!5] (10.8,5) \myupperbidiagonal;
      \def\nits{10}
      \def\nitsmone{9}
      \draw (10.8,5) rectangle ++(0.2*\nits,-0.2*\nits) rectangle +(0.2,0.2*\nits);
      \node[left] at (12.8,4.4) {$J_k^T$};
      \begin{scope}[xshift=-6cm,yshift=-2.8cm]
        \draw[fill=black!5] (16.8,5) \mylowerbidiagonalextra;
        \draw (16.8,5) rectangle +(0.2*\nits,-0.2*\nits);
        \node[left] at (18.8,4.4) {$\check{J}_k^T$};
      \end{scope}
      \node[below,font=\footnotesize] at (10,0) {(1)};
    \end{scope}
    \begin{scope}[xshift=5cm]
      \draw[fill=black!5] (8.3,0) rectangle node {$\widetilde{V}'_k$} +(2,5);
      \draw[fill=black!5] (10.3,0) rectangle +(0.2,5);
      \node[right,rotate=90] at (10.4,5) {$\tilde{v}_{k+1}$};
      \def\nits{11}
      \def\nitsmone{10}
      \draw[fill=black!5] (10.8,5) \mydiagonal;
      \def\nits{10}
      \def\nitsmone{9}
      \draw (10.8,5) rectangle ++(0.2*\nits,-0.2*\nits) rectangle +(0.2,0.2*\nits);
      \node[left] at (12.8,4.4) {$C_k$};
      \draw[fill=black!5] (10.8,2.8) rectangle +(2.2,0.2);
      \begin{scope}[xshift=-6cm,yshift=-2.8cm]
        \draw[fill=black!5] (16.8,5) \mydiagonal;
        \draw (16.8,5) rectangle +(0.2*\nits,-0.2*\nits);
        \node[left] at (18.8,4.4) {$S_k$};
        \draw[fill=black!5] (16.8,2.8) rectangle +(2,0.2);
      \end{scope}
      \node[below,font=\footnotesize] at (10,0) {(2)};
    \end{scope}
    \begin{scope}[xshift=10cm]
      \draw[fill=black!5] (8.3,0) rectangle node {$\tilde{V}'_r$} +(1.2,5);
      \draw[fill=black!5] (10.3,0) rectangle +(0.2,5);
      \draw[fill=black!25] (8.3,0) ++(1.2,0) rectangle +(0.8,5);
      \node[right,rotate=90] at (10.4,5) {$\tilde{v}_{k+1}$};
      \def\nits{11}
      \def\nitsmone{10}
      \draw[fill=black!5] (10.8,5) \mydiagonal;
      \def\nits{10}
      \def\nitsmone{9}
      \node[left] at (12.8,4.4) {$C_k$};
      \draw[fill=black!5] (10.8,2.8) rectangle +(2.2,0.2);
      \draw (10.8,5) rectangle ++(0.2*\nits,-0.2*\nits) rectangle +(0.2,0.2*\nits);
      \def\nits{4}
      \draw[fill=black!25] (12,3.8) \mydiagonal;
      \draw[fill=black!25] (12,2.8) rectangle +(0.8,0.2);
      \begin{scope}[xshift=-6cm,yshift=-2.8cm]
        \def\nits{10}
        \draw (16.8,5) rectangle +(0.2*\nits,-0.2*\nits);
        \draw[fill=black!5] (16.8,2.8) rectangle +(2,0.2);
        \draw[fill=black!5] (16.8,5) \mydiagonal;
        \draw[fill=black!25] (18,2.8) rectangle +(0.8,0.2);
        \node[left] at (18.8,4.4) {$S_k$};
        \def\nits{4}
        \draw[fill=black!25] (18,3.8) \mydiagonal;
      \end{scope}
      \node[below,font=\footnotesize] at (10,0) {(3)};
    \end{scope}
    \begin{scope}[xshift=15cm]
      \draw[fill=black!5] (8.3,0) rectangle node {$\tilde{V}'_r$} +(1.2,5);
      \draw[fill=black!5] (9.5,0) rectangle +(0.2,5);
      \node[right,rotate=90] at (9.7,5) {$\tilde{v}_{k+1}$};
      \def\nits{6}
      \draw[fill=black!5] (10,5) \mydiagonal;
      \node at (11.1,4.6) {$C_r$};
      \draw[fill=black!5] (10,3.6) rectangle ++(1.2,0.2) rectangle +(0.2,-0.2);
      \begin{scope}[xshift=-6cm,yshift=-2.8cm]
        \draw[fill=black!5] (16,5) \mydiagonal;
        \node at (17.1,4.6) {$S_r$};
        \draw[fill=black!5] (16,3.6) rectangle +(1.2,0.2);
      \end{scope}
      \node[below,font=\footnotesize] at (9.5,0) {(4)};
    \end{scope}
    \begin{scope}[xshift=18.5cm]
      \draw[fill=black!5] (8.3,0) rectangle node {$\tilde{V}'_r$} +(1.2,5);
      \draw[fill=black!5] (9.5,0) rectangle +(0.2,5);
      \node[right,rotate=90] at (9.7,5) {$\tilde{v}_{k+1}$};
      \draw[fill=black!25] (8.3,0) ++(1.4,0) rectangle +(0.8,5);
      \def\nits{6}
      \draw[fill=black!5] (10.8,5) \mydiagonal;
      \draw[fill=black!5] (10.8,3.6) rectangle ++(1.2,0.2) rectangle +(0.2,-0.2);
      \def\nits{5}
      \def\nitsmone{4}
      \draw[fill=black!25] (12,3.8) \myupperbidiagonal;
      \draw (10.8,5) rectangle ++(2,-2) rectangle +(0.2,2);
      \draw[fill=black!5] (12,3.8) rectangle +(0.2,-0.2);
      \begin{scope}[xshift=-6cm,yshift=-2.8cm]
        \def\nits{6}
        \draw[fill=black!5] (16.8,5) \mydiagonal;
        \draw[fill=black!5] (16.8,3.6) rectangle +(1.2,0.2);
        \def\nits{4}
        \draw[fill=black!25] (18,3.8) \mylowerbidiagonalextra;
        \draw (16.8,5) rectangle +(2,-2);
      \end{scope}
      \node[below,font=\footnotesize] at (10,0) {(5)};
    \end{scope}
  \end{tikzpicture}}
  \caption{Illustration of the five steps of thick restart in GSVD: (1) initial Lanczos factorization of order $k$, (2) solve projected problem, (3) sort and check convergence, (4) truncate to factorization of order $r$, and (5) extend to a factorization of order $k$.}
  \label{fig:thick-gsvd}
\end{figure}
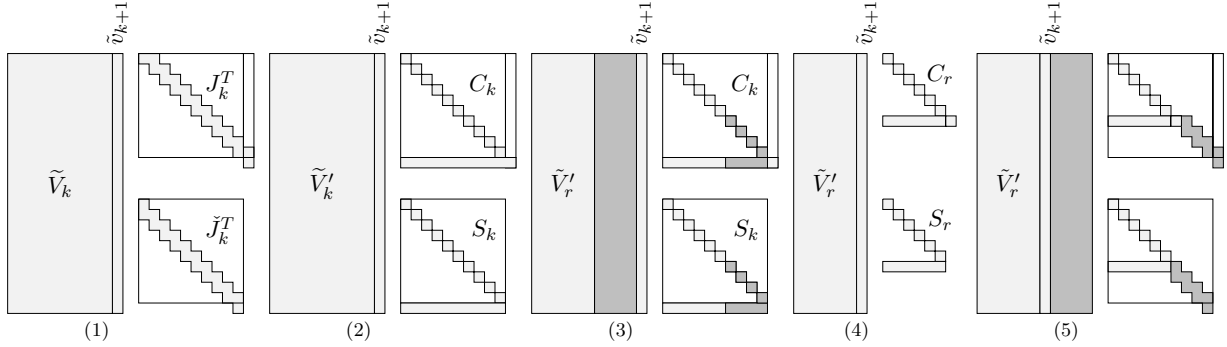

The process is illustrated graphically in \cref{fig:thick-gsvd}, in a similar way as was done for the SVD in \cref{fig:thick-svd}. This time, the pictures show the right Lanczos basis $\widetilde{V}_k$, together with both bidiagonals $J_k$ and $\check{J}_k$, and depict how these matrices evolve when the compression and extension are carried out. In the last panel (step 5), we can see how after restart both the upper and lower bidiagonals have the leading part with an arrowhead shape, with the spike pointing to the left in both cases.

\refstepcounter{myalg}
\begin{algorithm}
  \caption{Thick-restarted joint Lanczos bidiagonalization for the GSVD}
  \label{alg:restart}
  \begin{algorithmic}[1]
    \REQUIRE Matrices $A\in\mathbb{R}^{m\times n}$, $B\in\mathbb{R}^{p\times n}$, unit-norm vector $u_1\in\mathbb{R}^m$, maximum basis size $k$, restart size $r$, number of wanted generalized singular quadruples $s$, tolerance \texttt{tol}.
    \ENSURE Partial GSVD ${U'}_s^TAG_s=C_s,\; \widehat{U}_s^{\prime T}BG_s=S_s$.
    \STATE Compute $\hat{\beta}_0,\tilde{v}_1,\alpha_1$ as in \cref{alg:kilmer}
    \STATE $j'=0$
    \LOOP
    \FOR{$j=j'+1,j'+2,\dots,k$}
    \STATE Do step $j$ of bidiagonalization (Lines \ref{alg:kilmer:matlab1}--\ref{alg:kilmer:last-in-for} of \cref{alg:kilmer})
    \ENDFOR
    \STATE Compute the CS decomposition of $\lbrace J_k,\check{J}_k \rbrace$ as in \eqref{eq-J-CS-decomp}
    \STATE Reorder the CS decomposition according to $c_i/s_i$ (in increasing or decreasing order)
    \STATE Declare convergence if $\sqrt{\left(\alpha_{k+1} e_{k+1}^T x_i \right)^2 + \left(\check{\beta}_{k} e_{k}^T \hat{x}_i \right)^2}<\texttt{tol}$ for $i=1,\dots,s$, see \cref{sec:estim}
    \STATE Compute $U'_{r+1}$, $\widehat{U}'_{r}$, $\widetilde{V}'_{r}$, $b_{r+1}$ and $\hat{b}_r$ as in~\eqref{eq-restart-a}-\eqref{eq-restart-b}
    \STATE Exit loop if converged, setting $g_j,\;j=1,\dots,s,$ as the solution to the least squares problem 
    $\begin{bmatrix}A\\B\end{bmatrix}g_j=\tilde{v}'_j$
    \STATE $j'=r$
    \ENDLOOP
  \end{algorithmic}
\end{algorithm}

\Cref{alg:restart} summarizes the overall restarted procedure.

\subsection{Residual estimates}\label{sec:estim}

We now discuss a convergence criterion to be used in the computation of the partial GSVD with the method of \cref{sec:trkilmer} to decide when an approximate generalized singular quadruple $(\tilde{\sigma}_i,\tilde{u}^A_i,\tilde{u}^B_i,\tilde{g}_i)$, with $\tilde{\sigma}_i=\tilde{c}_i/\tilde{s}_i$, can be considered converged. We must derive formulas for the estimates of residual norms, using the information obtained during the bidiagonalization, without expensive additional computation. The iteration stops when the number of converged solutions reaches the number requested by the user.

As in the case of the SVD (see \cref{sec:residual}), there are two possible alternatives to define the residual associated to an approximate solution, in terms of either the cross product eigenproblem~\eqref{eq:gsvdeigcross} or the cyclic eigenproblem~\eqref{eq:gsvdeigcyclic},
\begin{align}
r_i^\mathrm{CROSS}&=\tilde{s}_i^2A^TA\tilde{g}_i-\tilde{c}_i^2B^TB\tilde{g}_i,
\label{eq:residualgsvdcross}
\\
r_i^\mathrm{CYCLIC,A}&=\begin{bmatrix}A\tilde{g}_i/\tilde{s}_i-\tilde{\sigma}_i\tilde{u}^A_i\\A^T\tilde{u}^A_i-\tilde{\sigma}_iB^TB\tilde{g}_i/\tilde{s}_i\end{bmatrix},
\label{eq:residualgsvdcyclic}
\\
r_i^\mathrm{CYCLIC,B}&=\begin{bmatrix}B\tilde{g}_i/\tilde{c}_i-\tilde{\sigma}_i^{-1}\tilde{u}^B_i\\B^T\tilde{u}^B_i-\tilde{\sigma}_i^{-1}A^TA\tilde{g}_i/\tilde{c}_i\end{bmatrix}.
\label{eq:residualgsvdcyclicb}
\end{align}

The residual norm $\left\|r_i^\mathrm{CROSS}\right\|_2$ was used by Zha~\cite{Zha:1996:CGS}, who showed, in the context of joint bidiagonalization with both bidiagonal matrices in upper form, that
$$\left\|(\tilde{s}_i^2 A^T A - \tilde{c}_i^2 B^T B) \tilde{g}_i\right\|_2 \le \alpha_k \beta_k \left|e_k^T y_i\right| \left\|R\right\|_2,$$
where $y_i$ is the $i$th right singular vector of the bidiagonal matrix $J_k$ resulting from a $k$-step joint bidiagonalization process.
In the case of joint bidiagonalization with lower-upper bidiagonal forms, the above inequality should be modified slightly,
\begin{equation}
\left\|(\tilde{s}_i^2 A^T A - \tilde{c}_i^2 B^T B) \tilde{g}_i\right\|_2 \le \alpha_{k+1} \beta_{k+1} \left|e_k^T y_i\right| \left\|R\right\|_2.
\label{eq-rnorm-bound-zha}
\end{equation}
The matrix $R$ above is the triangular factor of the QR decomposition~\eqref{eq-QR}, so its 2-norm is not readily available, but it can be bounded by $\sqrt{m+p}\max\{\|A\|_\infty,\|B\|_\infty\}$.

The residual $r_i^\mathrm{CROSS}$~\eqref{eq:residualgsvdcross} does not take into account possible errors in left vectors $\tilde{u}_i^A$, $\tilde{u}_i^B$, so the residual $r_i^\mathrm{CYCLIC,A}$~\eqref{eq:residualgsvdcyclic} may be more appropriate. This is the residual employed in the Jacobi-Davidson method~\cite{Hochstenbach:2009:JTM}. It refers to $\tilde{u}_i^A$. Similarly, the residual $r_i^\mathrm{CYCLIC,B}$~\eqref{eq:residualgsvdcyclicb} stemming from the second cyclic eigenproblem in~\eqref{eq:gsvdeigcyclic} refers to $\tilde{u}_i^B$. In the following, we define a convergence criterion that combines the two cyclic residuals so that both $\tilde{u}_i^A$ and $\tilde{u}_i^B$ are taken into consideration.

From the definition of the GSVD~\eqref{eq-gsvd} we have
\begin{align}
A g_i &= c_i u_i^A,              \label{eq-gsvd-Agu}\\
B g_i &= s_i u_i^B,              \label{eq-gsvd-Bgu}\\
s_i  A^T u_i^A &= c_i B^T u_i^B. \label{eq-gsvd-AtBt}
\end{align}
Using \eqref{eq-gsvd-Bgu} in \eqref{eq-gsvd-AtBt}, and in the same way using \eqref{eq-gsvd-Agu} in \eqref{eq-gsvd-AtBt}, we get
\begin{align}
s_i^2 A^T u_i^A = c_i B^T B g_i,\label{eq-gsvd-res1}\\
c_i^2 B^T u_i^B = s_i A^T A g_i,\label{eq-gsvd-res2}
\end{align}
and we define a residual that accounts for these two relations, whose norm is
\begin{equation}\label{eq:residual}
\left\|r^\mathrm{GSVD}\right\|_2=\sqrt{\left\|\tilde{s}_i^2 A^T \tilde{u}_i^A - \tilde{c}_i B^T B \tilde{g}_i\right\|_2^2 +
\left\|\tilde{c}_i^2 B^T \tilde{u}_i^B - \tilde{s}_i A^T A \tilde{g}_i\right\|_2^2}.
\end{equation}

We can derive bounds for this residual from the quantities computed in the joint bidiagonalization.
From \eqref{eq-kilmer-vectors-QA}-\eqref{eq-kilmer-vectors-QB} we have
\begin{align}
A \tilde{g}_i &= \tilde{c}_i \tilde{u}_i^A, &
A^T \tilde{u}_i^A &= R^T \left(\tilde{c}_i R \tilde{g}_i + \alpha_{k+1} v_{k+1} e_{k+1}^T x_i\right),
\label{eq-kilmer-final-vectors-A}
\\
B \tilde{g}_i &= \tilde{s}_i \tilde{u}_i^B, &
B^T \tilde{u}_i^B &= R^T \left(\tilde{s}_i R \tilde{g}_i + \check{\beta}_{k} v_{k+1} e_{k}^T \hat{x}_i\right).
\label{eq-kilmer-final-vectors-B}
\end{align}
Using~\eqref{eq-kilmer-final-vectors-A} and taking into account that $\tilde{c}_i^2+\tilde{s}_i^2=1$, the residual norm for~\eqref{eq-gsvd-res1} is
\begin{align*}
\left\|\tilde{s}_i^2 A^T \tilde{u}_i^A - \tilde{c}_i B^T B \tilde{g}_i\right\|_2 &=
\left\|A^T \tilde{u}_i^A - \tilde{c}_i Z^T Z \tilde{g}_i\right\|_2 \\ &=
\left\|A^T \tilde{u}_i^A - \tilde{c}_i R^T R \tilde{g}_i\right\|_2 = 
\left\|\alpha_{k+1} R^T v_{k+1} e_{k+1}^T x_i\right\|_2
\end{align*}
with the bound
\begin{align*}
\left\|\tilde{s}_i^2 A^T \tilde{u}_i^A - \tilde{c}_i B^T B \tilde{g}_i\right\|_2 \le
\alpha_{k+1} \left|e_{k+1}^T x_i \right| \left\|R \right\|_2 .
\end{align*}
In fact, this is a bound for $\tilde{s}_i^2\|r_i^\mathrm{CYCLIC,A}\|$ for the cyclic residual~\eqref{eq:residualgsvdcyclic}, since the upper part of that residual is exactly zero due to the left relation in~\eqref{eq-kilmer-final-vectors-A}.
Analogously, the residual associated to~\eqref{eq-gsvd-res2} verifies
\begin{align*}
\left\|\tilde{c}_i^2 B^T \tilde{u}_i^B - \tilde{s}_i A^T A \tilde{g}_i\right\|_2 \le
\left|\check{\beta}_{k} e_{k}^T \hat{x}_i \right| \left\|R \right\|_2 ,
\end{align*}
which is related to $\tilde{c}_i^2\|r_i^\mathrm{CYCLIC,B}\|$. Combining the previous two bounds,
\begin{equation}
\left\|r^\mathrm{GSVD}\right\|_2 \le \sqrt{\left(\alpha_{k+1} e_{k+1}^T x_i \right)^2 + \left(\check{\beta}_{k} e_{k}^T \hat{x}_i \right)^2} \left\|R \right\|_2.
\label{eq-rgsvd-bound}
\end{equation}

It is interesting to see that $\left\|r^\mathrm{GSVD}\right\|_2$ is equal to the residual norm of~\eqref{eq-gsvd-AtBt}. Using~\eqref{eq-gsvd-Agu} and~\eqref{eq-gsvd-Bgu} we have
\begin{align*}
\left\|\tilde{s}_i^2 A^T \tilde{u}_i^A - \tilde{c}_i B^T B \tilde{g}_i\right\|_2 = \left\|\tilde{s}_i^2 A^T \tilde{u}_i^A - \tilde{c}_i \tilde{s}_i B^T \tilde{u}_i^B \right\|_2 = \tilde{s}_i \left\|\tilde{s}_i A^T \tilde{u}_i^A - \tilde{c}_i B^T \tilde{u}_i^B \right\|_2,\\
\left\|\tilde{c}_i^2 B^T \tilde{u}_i^B - \tilde{s}_i A^T A \tilde{g}_i\right\|_2 = \left\|\tilde{c}_i^2 B^T \tilde{u}_i^B - \tilde{s}_i \tilde{c}_i A^T \tilde{u}_i^A \right\|_2 = \tilde{c}_i \left\|\tilde{c}_i B^T \tilde{u}_i^B - \tilde{s}_i A^T \tilde{u}_i^A \right\|_2,
\end{align*}
and
$$
\left\|r^\mathrm{GSVD}\right\|_2 = \left\|\tilde{s}_i A^T \tilde{u}_i^A - \tilde{c}_i B^T \tilde{u}_i^B \right\|_2.
$$

As a summary, our criterion to accept a computed generalized singular quadruple as converged during the restart of the joint bidiagonalization is $\sqrt{\left(\alpha_{k+1} e_{k+1}^T x_i \right)^2 + \left(\check{\beta}_{k} e_{k}^T \hat{x}_i \right)^2}<\texttt{tol}$, where \texttt{tol} is the user-defined tolerance. According to~\eqref{eq-rgsvd-bound}, this can be considered a criterion relative to the norm of the problem matrices.

\subsection{Using a scale factor}\label{sec:scaling}

We have implemented an option for scaling one of the matrices, as this may be beneficial for convergence in some cases. When the user specifies a scale factor $\gamma$, the method is applied to the pair $\{A,\gamma B\}$ instead of $\{A,B\}$. As a consequence, the algorithm works with a modified coefficient matrix for the least squares problem,
\begin{equation*}
Z_\gamma:=
\begin{bmatrix}
A\\
\gamma B
\end{bmatrix}
=
\begin{bmatrix}
Q_{A,\gamma}\\
Q_{B,\gamma}
\end{bmatrix}
R_\gamma,
\label{eq-QR-scaled}
\end{equation*}
where $Q_{A,\gamma}$ and $Q_{B,\gamma}$ are different from those obtained without scaling.

If we use $(c_i,s_i,u^A_i,u^B_i,g_i)$ to denote the solution for the pair $\{A,B\}$ as in~\eqref{eq:gsvd2}, we can write the relations for the scaled problem as
\begin{equation}
\label{eq:gsvd2-scaled}
A (\omega_i g_i)=(\omega_i c_i) u_i^A,\qquad \gamma B (\omega_i g_i)=(\gamma \omega_i s_i) u_i^B,
\end{equation}
for certain weights $\omega_i$, so after solving~\eqref{eq:gsvd2-scaled} we can retrieve the solution of the original problem~\eqref{eq:gsvd2} by multiplying the generalized singular values $\sigma_i$ by $\gamma$ and also multiplying the $g_i$ vectors by $\omega_i^{-1}$, which can be obtained from the relation
\begin{equation*}
\omega_i^2c_i^2+\gamma^2\omega_i^2s_i^2=1,\quad
\omega_i^{-1}=\sqrt{c_i^2+\gamma^2s_i^2},
\end{equation*}
where $c_i, s_i$, corresponding to the original problem, can be computed from $\sigma_i$.

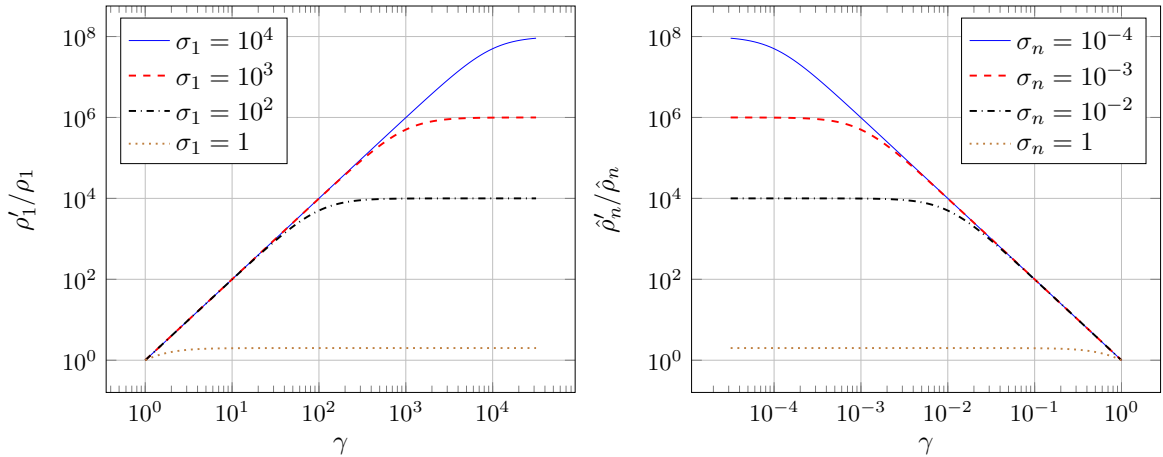
\begin{figure}
	\centering
   \begin{tikzpicture}[scale=0.9]
		\begin{loglogaxis}[
		xlabel={$\gamma$},
		ylabel={$\rho'_1/\rho_1$},
		grid=major,no markers,
		ticklabel style={font=\small},
		legend style={anchor=north, cells={anchor=west}}, legend pos=north west,
		cycle list={
		   {solid,blue},
		   {dashed,red,thick},
		   {dashdotted,black,thick},
		   {dotted,brown,thick},
		},
		]
		\pgfplotstableread{
   1.0000000e+00   1.0000000e+00   1.0000000e+00   1.0000000e+00   1.0000000e+00
   1.1103363e+00   1.1042822e+00   1.2328178e+00   1.2328462e+00   1.2328465e+00
   1.2328467e+00   1.2063208e+00   1.5198316e+00   1.5199098e+00   1.5199106e+00
   1.3688745e+00   1.3040610e+00   1.8736528e+00   1.8738149e+00   1.8738165e+00
   1.5199111e+00   1.3957933e+00   2.3098258e+00   2.3101254e+00   2.3101284e+00
   1.6876125e+00   1.4802534e+00   2.8475078e+00   2.8480288e+00   2.8480340e+00
   1.8738174e+00   1.5566571e+00   3.5103078e+00   3.5111804e+00   3.5111891e+00
   2.0805675e+00   1.6246770e+00   4.3273176e+00   4.3287435e+00   4.3287578e+00
   2.3101297e+00   1.6843769e+00   5.3343818e+00   5.3366718e+00   5.3366947e+00
   2.5650209e+00   1.7361230e+00   6.5756583e+00   6.5792900e+00   6.5793263e+00
   2.8480359e+00   1.7804910e+00   8.1055377e+00   8.1112435e+00   8.1113006e+00
   3.1622777e+00   1.8181802e+00   9.9910000e+00   9.9999010e+00   9.9999901e+00
   3.5111917e+00   1.8499435e+00   1.2314507e+01   1.2328316e+01   1.2328455e+01
   3.8986037e+00   1.8765347e+00   1.5177548e+01   1.5198881e+01   1.5199094e+01
   4.3287613e+00   1.8986717e+00   1.8704981e+01   1.8737824e+01   1.8738153e+01
   4.8063809e+00   1.9170151e+00   2.3050336e+01   2.3100764e+01   2.3101270e+01
   5.3366992e+00   1.9321564e+00   2.8402289e+01   2.8479549e+01   2.8480323e+01
   5.9255310e+00   1.9446147e+00   3.4992529e+01   3.5110686e+01   3.5111871e+01
   6.5793322e+00   1.9548387e+00   4.3105307e+01   4.3285740e+01   4.3287552e+01
   7.3052715e+00   1.9632110e+00   5.3088957e+01   5.3364145e+01   5.3366912e+01
   8.1113083e+00   1.9700549e+00   6.5369748e+01   6.5788995e+01   6.5793215e+01
   9.0062802e+00   1.9756414e+00   8.0468410e+01   8.1106505e+01   8.1112938e+01
   1.0000000e+01   1.9801961e+00   9.9019704e+01   9.9990002e+01   9.9999802e+01
   1.1103363e+01   1.9839059e+00   1.2179533e+02   1.2326948e+02   1.2328440e+02
   1.2328467e+01   1.9869254e+00   1.4973039e+02   1.5196801e+02   1.5199073e+02
   1.3688745e+01   1.9893813e+00   1.8395334e+02   1.8734664e+02   1.8738121e+02
   1.5199111e+01   1.9913778e+00   2.2581913e+02   2.3095962e+02   2.3101221e+02
   1.6876125e+01   1.9930002e+00   2.7694431e+02   2.8472250e+02   2.8480249e+02
   1.8738174e+01   1.9943181e+00   3.3924248e+02   3.5099593e+02   3.5111759e+02
   2.0805675e+01   1.9953884e+00   4.1495651e+02   4.3268883e+02   4.3287383e+02
   2.3101297e+01   1.9962574e+00   5.0668263e+02   5.3338527e+02   5.3366655e+02
   2.5650209e+01   1.9969628e+00   6.1737895e+02   6.5750063e+02   6.5792825e+02
   2.8480359e+01   1.9975353e+00   7.5034809e+02   8.1047343e+02   8.1112345e+02
   3.1622777e+01   1.9980000e+00   9.0918091e+02   9.9900100e+02   9.9998901e+02
   3.5111917e+01   1.9983771e+00   1.0976459e+03   1.2313287e+03   1.2328303e+03
   3.8986037e+01   1.9986830e+00   1.3195081e+03   1.5176045e+03   1.5198865e+03
   4.3287613e+01   1.9989312e+00   1.5782649e+03   1.8703128e+03   1.8737805e+03
   4.8063809e+01   1.9991326e+00   1.8767945e+03   2.3048053e+03   2.3100740e+03
   5.3366992e+01   1.9992960e+00   2.2169286e+03   2.8399476e+03   2.8479519e+03
   5.9255310e+01   1.9994286e+00   2.5989856e+03   3.4989064e+03   3.5110650e+03
   6.5793322e+01   1.9995361e+00   3.0213287e+03   4.3101039e+03   4.3285696e+03
   7.3052715e+01   1.9996233e+00   3.4800367e+03   5.3083701e+03   5.3364092e+03
   8.1113083e+01   1.9996941e+00   3.9687869e+03   6.5363276e+03   6.5788929e+03
   9.0062802e+01   1.9997515e+00   4.4790311e+03   8.0460444e+03   8.1106424e+03
   1.0000000e+02   1.9997980e+00   5.0004950e+03   9.9009901e+03   9.9989902e+03
   1.1103363e+02   1.9998358e+00   5.5219589e+03   1.2178327e+04   1.2326935e+04
   1.2328467e+02   1.9998664e+00   6.0322031e+03   1.4971556e+04   1.5196786e+04
   1.3688745e+02   1.9998913e+00   6.5209533e+03   1.8393513e+04   1.8734645e+04
   1.5199111e+02   1.9999114e+00   6.9796613e+03   2.2579677e+04   2.3095939e+04
   1.6876125e+02   1.9999278e+00   7.4020044e+03   2.7691689e+04   2.8472221e+04
   1.8738174e+02   1.9999410e+00   7.7840614e+03   3.3920890e+04   3.5099558e+04
   2.0805675e+02   1.9999518e+00   8.1241955e+03   4.1491543e+04   4.3268840e+04
   2.3101297e+02   1.9999605e+00   8.4227251e+03   5.0663247e+04   5.3338474e+04
   2.5650209e+02   1.9999676e+00   8.6814819e+03   6.1731783e+04   6.5749998e+04
   2.8480359e+02   1.9999733e+00   8.9033441e+03   7.5027381e+04   8.1047263e+04
   3.1622777e+02   1.9999780e+00   9.0918091e+03   9.0909091e+04   9.9900001e+04
   3.5111917e+02   1.9999818e+00   9.2506419e+03   1.0975372e+05   1.2313275e+05
   3.8986037e+02   1.9999848e+00   9.3836111e+03   1.3193774e+05   1.5176030e+05
   4.3287613e+02   1.9999873e+00   9.4943074e+03   1.5781087e+05   1.8703109e+05
   4.8063809e+02   1.9999893e+00   9.5860335e+03   1.8766087e+05   2.3048030e+05
   5.3366992e+02   1.9999910e+00   9.6617475e+03   2.2167092e+05   2.8399448e+05
   5.9255310e+02   1.9999923e+00   9.7240457e+03   2.5987284e+05   3.4989029e+05
   6.5793322e+02   1.9999934e+00   9.7751709e+03   3.0210297e+05   4.3100996e+05
   7.3052715e+02   1.9999943e+00   9.8170367e+03   3.4796922e+05   5.3083648e+05
   8.1113083e+02   1.9999950e+00   9.8512596e+03   3.9683940e+05   6.5363211e+05
   9.0062802e+02   1.9999955e+00   9.8791947e+03   4.4785877e+05   8.0460364e+05
   1.0000000e+03   1.9999960e+00   9.9019703e+03   5.0000000e+05   9.9009803e+05
   1.1103363e+03   1.9999964e+00   9.9205216e+03   5.5214123e+05   1.2178315e+06
   1.2328467e+03   1.9999967e+00   9.9356203e+03   6.0316060e+05   1.4971542e+06
   1.3688745e+03   1.9999969e+00   9.9479010e+03   6.5203078e+05   1.8393495e+06
   1.5199111e+03   1.9999971e+00   9.9578847e+03   6.9789703e+05   2.2579655e+06
   1.6876125e+03   1.9999973e+00   9.9659975e+03   7.4012716e+05   2.7691662e+06
   1.8738174e+03   1.9999974e+00   9.9725877e+03   7.7832908e+05   3.3920856e+06
   2.0805675e+03   1.9999975e+00   9.9779397e+03   8.1233913e+05   4.1491502e+06
   2.3101297e+03   1.9999976e+00   9.9822850e+03   8.4218913e+05   5.0663197e+06
   2.5650209e+03   1.9999977e+00   9.9858125e+03   8.6806226e+05   6.1731722e+06
   2.8480359e+03   1.9999978e+00   9.9886755e+03   8.9024628e+05   7.5027306e+06
   3.1622777e+03   1.9999978e+00   9.9909990e+03   9.0909091e+05   9.0909001e+06
   3.5111917e+03   1.9999978e+00   9.9928845e+03   9.2497262e+05   1.0975361e+07
   3.8986037e+03   1.9999979e+00   9.9944143e+03   9.3826822e+05   1.3193761e+07
   4.3287613e+03   1.9999979e+00   9.9956556e+03   9.4933675e+05   1.5781071e+07
   4.8063809e+03   1.9999979e+00   9.9966627e+03   9.5850846e+05   1.8766069e+07
   5.3366992e+03   1.9999979e+00   9.9974797e+03   9.6607911e+05   2.2167069e+07
   5.9255310e+03   1.9999979e+00   9.9981425e+03   9.7230831e+05   2.5987258e+07
   6.5793322e+03   1.9999980e+00   9.9986802e+03   9.7742032e+05   3.0210266e+07
   7.3052715e+03   1.9999980e+00   9.9991163e+03   9.8160649e+05   3.4796887e+07
   8.1113083e+03   1.9999980e+00   9.9994702e+03   9.8502844e+05   3.9683900e+07
   9.0062802e+03   1.9999980e+00   9.9997572e+03   9.8782167e+05   4.4785832e+07
   1.0000000e+04   1.9999980e+00   9.9999900e+03   9.9009901e+05   4.9999950e+07
   1.1103363e+04   1.9999980e+00   1.0000179e+04   9.9195396e+05   5.5214067e+07
   1.2328467e+04   1.9999980e+00   1.0000332e+04   9.9346367e+05   6.0315999e+07
   1.3688745e+04   1.9999980e+00   1.0000456e+04   9.9469163e+05   6.5203013e+07
   1.5199111e+04   1.9999980e+00   1.0000557e+04   9.9568990e+05   6.9789633e+07
   1.6876125e+04   1.9999980e+00   1.0000639e+04   9.9650109e+05   7.4012642e+07
   1.8738174e+04   1.9999980e+00   1.0000705e+04   9.9716005e+05   7.7832830e+07
   2.0805675e+04   1.9999980e+00   1.0000759e+04   9.9769519e+05   8.1233830e+07
   2.3101297e+04   1.9999980e+00   1.0000803e+04   9.9812969e+05   8.4218828e+07
   2.5650209e+04   1.9999980e+00   1.0000838e+04   9.9848240e+05   8.6806138e+07
   2.8480359e+04   1.9999980e+00   1.0000867e+04   9.9876867e+05   8.9024538e+07
   3.1622777e+04   1.9999980e+00   1.0000890e+04   9.9900100e+05   9.0908999e+07
   }\scalingdata
		\addplot table [y index=4]{\scalingdata};
		\addlegendentry{$\sigma_1=10^4$}
		\addplot table [y index=3]{\scalingdata};
		\addlegendentry{$\sigma_1=10^3$}
		\addplot table [y index=2]{\scalingdata};
		\addlegendentry{$\sigma_1=10^2$}
		\addplot table [y index=1]{\scalingdata};
		\addlegendentry{$\sigma_1=1$}
		\end{loglogaxis}
	\end{tikzpicture}
   \begin{tikzpicture}[scale=0.9]
		\begin{loglogaxis}[
		xlabel={$\gamma$},
		ylabel={$\hat{\rho}'_n/\hat{\rho}_n$},
		grid=major,no markers,
		ticklabel style={font=\small},
		legend style={anchor=north, cells={anchor=west}}, legend pos=north east,
		cycle list={
		   {solid,blue},
		   {dashed,red,thick},
		   {dashdotted,black,thick},
		   {dotted,brown,thick},
		},
		]
		\pgfplotstableread{
   3.1622777e-05   1.9998000e+00   9.9999000e+03   9.9890211e+05   9.0900002e+07
   3.5111917e-05   1.9998000e+00   9.9998767e+03   9.9866980e+05   8.9015727e+07
   3.8986037e-05   1.9998000e+00   9.9998480e+03   9.9838356e+05   8.6797547e+07
   4.3287613e-05   1.9998000e+00   9.9998126e+03   9.9803088e+05   8.4210493e+07
   4.8063809e-05   1.9998000e+00   9.9997690e+03   9.9759643e+05   8.1225791e+07
   5.3366992e-05   1.9998000e+00   9.9997152e+03   9.9706134e+05   7.7825127e+07
   5.9255310e-05   1.9998000e+00   9.9996489e+03   9.9640245e+05   7.4005316e+07
   6.5793322e-05   1.9998000e+00   9.9995671e+03   9.9559133e+05   6.9782726e+07
   7.3052715e-05   1.9998000e+00   9.9994664e+03   9.9459317e+05   6.5196559e+07
   8.1113083e-05   1.9998000e+00   9.9993421e+03   9.9336533e+05   6.0310030e+07
   9.0062802e-05   1.9998000e+00   9.9991889e+03   9.9185576e+05   5.5208602e+07
   1.0000000e-04   1.9998000e+00   9.9990001e+03   9.9000100e+05   4.9995001e+07
   1.1103363e-04   1.9998000e+00   9.9987673e+03   9.8772389e+05   4.4781400e+07
   1.2328467e-04   1.9998000e+00   9.9984803e+03   9.8493094e+05   3.9679972e+07
   1.3688745e-04   1.9998000e+00   9.9981265e+03   9.8150932e+05   3.4793443e+07
   1.5199111e-04   1.9998000e+00   9.9976904e+03   9.7732357e+05   3.0207276e+07
   1.6876125e-04   1.9998000e+00   9.9971528e+03   9.7221206e+05   2.5984686e+07
   1.8738174e-04   1.9997999e+00   9.9964900e+03   9.6598348e+05   2.2164875e+07
   2.0805675e-04   1.9997999e+00   9.9956731e+03   9.5841357e+05   1.8764211e+07
   2.3101297e-04   1.9997999e+00   9.9946661e+03   9.4924278e+05   1.5779509e+07
   2.5650209e-04   1.9997999e+00   9.9934250e+03   9.3817534e+05   1.3192455e+07
   2.8480359e-04   1.9997999e+00   9.9918953e+03   9.2488106e+05   1.0974275e+07
   3.1622777e-04   1.9997998e+00   9.9900100e+03   9.0900092e+05   9.0900002e+06
   3.5111917e-04   1.9997998e+00   9.9876867e+03   8.9015815e+05   7.5019880e+06
   3.8986037e-04   1.9997997e+00   9.9848240e+03   8.6797633e+05   6.1725611e+06
   4.3287613e-04   1.9997996e+00   9.9812969e+03   8.4210576e+05   5.0658182e+06
   4.8063809e-04   1.9997996e+00   9.9769519e+03   8.1225871e+05   4.1487395e+06
   5.3366992e-04   1.9997995e+00   9.9716005e+03   7.7825204e+05   3.3917498e+06
   5.9255310e-04   1.9997993e+00   9.9650109e+03   7.4005390e+05   2.7688921e+06
   6.5793322e-04   1.9997992e+00   9.9568990e+03   6.9782795e+05   2.2577420e+06
   7.3052715e-04   1.9997990e+00   9.9469163e+03   6.5196624e+05   1.8391674e+06
   8.1113083e-04   1.9997987e+00   9.9346367e+03   6.0310089e+05   1.4970060e+06
   9.0062802e-04   1.9997984e+00   9.9195396e+03   5.5208657e+05   1.2177110e+06
   1.0000000e-03   1.9997980e+00   9.9009901e+03   4.9995050e+05   9.9000002e+05
   1.1103363e-03   1.9997976e+00   9.8782167e+03   4.4781444e+05   8.0452399e+05
   1.2328467e-03   1.9997970e+00   9.8502844e+03   3.9680012e+05   6.5356741e+05
   1.3688745e-03   1.9997963e+00   9.8160649e+03   3.4793477e+05   5.3078393e+05
   1.5199111e-03   1.9997954e+00   9.7742032e+03   3.0207306e+05   4.3096729e+05
   1.6876125e-03   1.9997943e+00   9.7230831e+03   2.5984711e+05   3.4985566e+05
   1.8738174e-03   1.9997930e+00   9.6607911e+03   2.2164897e+05   2.8396637e+05
   2.0805675e-03   1.9997914e+00   9.5850846e+03   1.8764230e+05   2.3045749e+05
   2.3101297e-03   1.9997893e+00   9.4933675e+03   1.5779525e+05   1.8701258e+05
   2.5650209e-03   1.9997869e+00   9.3826822e+03   1.3192468e+05   1.5174527e+05
   2.8480359e-03   1.9997838e+00   9.2497262e+03   1.0974286e+05   1.2312056e+05
   3.1622777e-03   1.9997800e+00   9.0909091e+03   9.0900092e+04   9.9890112e+04
   3.5111917e-03   1.9997754e+00   8.9024628e+03   7.5019954e+04   8.1039240e+04
   3.8986037e-03   1.9997696e+00   8.6806226e+03   6.1725673e+04   6.5743490e+04
   4.3287613e-03   1.9997626e+00   8.4218913e+03   5.0658232e+04   5.3333194e+04
   4.8063809e-03   1.9997538e+00   8.1233913e+03   4.1487436e+04   4.3264557e+04
   5.3366992e-03   1.9997431e+00   7.7832909e+03   3.3917532e+04   3.5096084e+04
   5.9255310e-03   1.9997298e+00   7.4012717e+03   2.7688948e+04   2.8469403e+04
   6.5793322e-03   1.9997135e+00   6.9789704e+03   2.2577442e+04   2.3093652e+04
   7.3052715e-03   1.9996933e+00   6.5203079e+03   1.8391693e+04   1.8732791e+04
   8.1113083e-03   1.9996685e+00   6.0316060e+03   1.4970075e+04   1.5195282e+04
   9.0062802e-03   1.9996378e+00   5.5214123e+03   1.2177122e+04   1.2325715e+04
   1.0000000e-02   1.9996001e+00   5.0000000e+03   9.9000101e+03   9.9980005e+03
   1.1103363e-02   1.9995535e+00   4.4785878e+03   8.0452480e+03   8.1098396e+03
   1.2328467e-02   1.9994961e+00   3.9683941e+03   6.5356806e+03   6.5782417e+03
   1.3688745e-02   1.9994254e+00   3.4796922e+03   5.3078447e+03   5.3358810e+03
   1.5199111e-02   1.9993382e+00   3.0210297e+03   4.3096773e+03   4.3281412e+03
   1.6876125e-02   1.9992307e+00   2.5987284e+03   3.4985601e+03   3.5107175e+03
   1.8738174e-02   1.9990982e+00   2.2167092e+03   2.8396666e+03   2.8476701e+03
   2.0805675e-02   1.9989348e+00   1.8766088e+03   2.3045772e+03   2.3098455e+03
   2.3101297e-02   1.9987335e+00   1.5781088e+03   1.8701278e+03   1.8735951e+03
   2.5650209e-02   1.9984853e+00   1.3193775e+03   1.5174543e+03   1.5197361e+03
   2.8480359e-02   1.9981794e+00   1.0975373e+03   1.2312069e+03   1.2327084e+03
   3.1622777e-02   1.9978024e+00   9.0909100e+02   9.9890221e+02   9.9989012e+02
   3.5111917e-02   1.9973379e+00   7.5027390e+02   8.1039330e+02   8.1104326e+02
   3.8986037e-02   1.9967654e+00   6.1731793e+02   6.5743565e+02   6.5786322e+02
   4.3287613e-02   1.9960601e+00   5.0663257e+02   5.3333257e+02   5.3361382e+02
   4.8063809e-02   1.9951913e+00   4.1491553e+02   4.3264610e+02   4.3283108e+02
   5.3366992e-02   1.9941213e+00   3.3920900e+02   3.5096129e+02   3.5108294e+02
   5.9255310e-02   1.9928036e+00   2.7691699e+02   2.8469441e+02   2.8477440e+02
   6.5793322e-02   1.9911815e+00   2.2579687e+02   2.3093685e+02   2.3098944e+02
   7.3052715e-02   1.9891854e+00   1.8393523e+02   1.8732819e+02   1.8736276e+02
   8.1113083e-02   1.9867300e+00   1.4971566e+02   1.5195307e+02   1.5197578e+02
   9.0062802e-02   1.9837111e+00   1.2178337e+02   1.2325737e+02   1.2327230e+02
   1.0000000e-01   1.9800020e+00   9.9010000e+01   9.9980203e+01   9.9990002e+01
   1.1103363e-01   1.9754482e+00   8.0460543e+01   8.1098576e+01   8.1105008e+01
   1.2328467e-01   1.9698629e+00   6.5363375e+01   6.5782582e+01   6.5786801e+01
   1.3688745e-01   1.9630203e+00   5.3083800e+01   5.3358962e+01   5.3361728e+01
   1.5199111e-01   1.9546497e+00   4.3101138e+01   4.3281554e+01   4.3283366e+01
   1.6876125e-01   1.9444277e+00   3.4989164e+01   3.5107309e+01   3.5108495e+01
   1.8738174e-01   1.9319718e+00   2.8399576e+01   2.8476828e+01   2.8477603e+01
   2.0805675e-01   1.9168335e+00   2.3048153e+01   2.3098577e+01   2.3099082e+01
   2.3101297e-01   1.8984938e+00   1.8703228e+01   1.8736068e+01   1.8736397e+01
   2.5650209e-01   1.8763611e+00   1.5176144e+01   1.5197475e+01   1.5197689e+01
   2.8480359e-01   1.8497753e+00   1.2313387e+01   1.2327195e+01   1.2327333e+01
   3.1622777e-01   1.8180182e+00   9.9901099e+00   9.9990101e+00   9.9990992e+00
   3.5111917e-01   1.7803365e+00   8.1048342e+00   8.1105396e+00   8.1105967e+00
   3.8986037e-01   1.7359773e+00   6.5751063e+00   6.5787377e+00   6.5787740e+00
   4.3287613e-01   1.6842414e+00   5.3339527e+00   5.3362425e+00   5.3362654e+00
   4.8063809e-01   1.6245533e+00   4.3269882e+00   4.3284140e+00   4.3284283e+00
   5.3366992e-01   1.5565468e+00   3.5100593e+00   3.5109318e+00   3.5109406e+00
   5.9255310e-01   1.4801583e+00   2.8473249e+00   2.8478458e+00   2.8478510e+00
   6.5793322e-01   1.3957149e+00   2.3096961e+00   2.3099957e+00   2.3099987e+00
   7.3052715e-01   1.3040008e+00   1.8735664e+00   1.8737284e+00   1.8737300e+00
   8.1113083e-01   1.2062799e+00   1.5197801e+00   1.5198583e+00   1.5198591e+00
   9.0062802e-01   1.1042616e+00   1.2327948e+00   1.2328232e+00   1.2328235e+00
   1.0000000e+00   1.0000000e+00   1.0000000e+00   1.0000000e+00   1.0000000e+00
   		}\scalingdata
		\addplot table [y index=4]{\scalingdata};
		\addlegendentry{$\sigma_n=10^{-4}$}
		\addplot table [y index=3]{\scalingdata};
		\addlegendentry{$\sigma_n=10^{-3}$}
		\addplot table [y index=2]{\scalingdata};
		\addlegendentry{$\sigma_n=10^{-2}$}
		\addplot table [y index=1]{\scalingdata};
		\addlegendentry{$\sigma_n=1$}
		\end{loglogaxis}
	\end{tikzpicture}
\caption{\label{fig:scaling}Left: relative gap ratio $\rho'_1/\rho_1$ vs scale factor $\gamma$, for $\sigma_n=10^{-3}$ and different values of $\sigma_1$. Right: relative gap ratio $\hat{\rho}'_n/\hat{\rho}_n$ vs scale factor $\gamma$, for $\sigma_1=100$ and different values of $\sigma_n$.}
\end{figure}

The convergence of the joint Lanczos bidiagonalization for the GSVD is determined by the convergence of the Lanczos bidiagonalization of $Q_A$
to approximate the values $c_i$ (or equivalently, the bidiagonalization of $Q_B$ to approximate the values $s_i$). As \cite{Parlett:1980:SEP} shows, the relative gap of the first two eigenvalues, $\rho_1=(\lambda_1-\lambda_2)/(\lambda_2-\lambda_{n})$, with $\lambda_i=c_i^2$ in our case, is of key importance for the convergence rate of Lanczos to approximate $\lambda_1$, with convergence increasing as $\rho_1$ grows.
Thus, it is interesting to analyze how $\rho_1$ is affected by scaling. 
The relative gap of the scaled problem is $\rho'_1=({c'_1}^2-{c'_2}^2)/({c'_2}^2-{c'_n}^2)$, with $c'_i=\omega_ic_i$, according to \eqref{eq:gsvd2-scaled}. The ratio $\rho'_1/\rho_1$ is
$$\frac{\rho'_1}{\rho_1}=\frac{({c'_1}^2-{c'_2}^2)(c_2^2-c_n^2)}{({c'_2}^2-{c'_n}^2)(c_1^2-c_2^2)}=\frac{c_n^2+\gamma^2(1-c_n^2)}{c_1^2+\gamma^2(1-c_1^2)}=\frac{\omega_1^2}{\omega_n^2}.$$
As $\gamma$ grows, $\rho'_1/\rho_1$ will also grow, tending to $(1-c_n^2)/(1-c_1^2)$, and favoring convergence of the scaled problem. This can be seen in \cref{fig:scaling} (left), which shows the value of $\rho'_1/\rho_1$ vs $\gamma$ for $\sigma_n=c_n/s_n=10^{-3}$ and different values of $\sigma_1=c_1/s_1$. Clearly, problems with large $\sigma_1$ should benefit the most from scaling. The figure suggests that values of $\gamma$ near $\sigma_1$ should be good for convergence (but of course $\sigma_1$ is unknown a priori).

On the other hand, if we want to compute the smallest generalized singular value $\sigma_n=c_n/s_n$, the relative gap of interest is $\hat{\rho}_n=(c_{n-1}^2-c_n^2)/(c_1^2-c_{n-1}^2)$, with the corresponding ratio
$$\frac{\hat{\rho}'_n}{\hat{\rho}_n}= \frac{(c_{n-1}^{\prime\,2}-c_n^{\prime \,2})(c_1^2-c_{n-1}^2)} {(c_{1}^{\prime\,2}-c_{n-1}^{\prime\,2})(c_{n-1}^2-c_n^2)}=\frac{\omega_n^2}{\omega_1^2}.$$
In this case, the ratio will grow as $\gamma$ decreases, tending to $c_1^2/c_n^2$, as can be seen in \cref{fig:scaling} (right).

\Cref{sec:results} includes some results comparing different values of $\gamma$. A good choice of $\gamma$ is application-dependent. An in-depth theoretical analysis of the impact of scaling on performance is out of the scope of this paper, as it depends on many factors. On one hand, a suitable value of $\gamma$ can improve the relative separation of the target singular values of $Q_{A,\gamma}$ and $Q_{B,\gamma}$, but on the other hand this may bring an increase of the condition number of matrix $Z_\gamma$ with the consequent negative effects (especially if using an iterative least squares solver such as LSQR).

\subsection{One-sided orthogonalization}\label{sec:oneside}

When the thick restart technique was introduced for the SVD in \cref{sec:restart}, we pointed out that the loss of orthogonality among Lanczos vectors can be avoided by means of full orthogonalization, that is, enforcing the full orthogonality of both left and right Lanczos bases via explicit orthogonalization. In practice, \cref{alg:bidiagsvd:expanda,alg:bidiagsvd:expandat} of \cref{alg:bidiagsvd} are replaced with $p_j=\mathrm{orthog}(Aq_j,P_{j-1})$ and $q_{j+1}=\mathrm{orthog}(A^Tp_j,Q_j)$. However, Simon and Zha~\cite{Simon:2000:LMA} showed that it is not necessary to explicitly orthogonalize both bases, and it is enough with orthogonalizing one of them since the level of orthogonality of left and right Lanczos vectors go hand in hand. This technique is called one-sided orthogonalization, and is implemented in SLEPc's thick-restart Lanczos solver for the SVD~\cite{Hernandez:2008:REP}, reducing the cost of each iteration significantly.

The one-sided orthogonalization technique can also be applied to the GSVD case. In Kilmer's joint bidiagonalization (\cref{alg:kilmer}), \cref{alg:kilmer:matlab1,alg:kilmer:matlab2,alg:kilmer:expand2} involve a full orthogonalization that would be implemented with a call to $\mathrm{orthog}()$. However, two of these three orthogonalizations can be avoided, that is, instead of explicitly orthogonalizing against the full basis, just a local orthogonalization with the previous vector is done. In particular, we have found that it is enough to explicitly orthogonalize the $U_{k+1}$ basis (\cref{alg:kilmer:matlab2}). This saves about two thirds of the orthogonalization cost, as will be illustrated in the computational results of \cref{sec:results}.

One-sided orthogonalization in the context of joint Lanczos bidiagonalization is discussed in~\cite{Jia:2023:JBM}, where the authors show that it is sufficient to explicitly orthogonalize $U_{k+1}$ and $V_k$, provided that the bidiagonal matrix $J_k$ does not become too ill-conditioned. Note that our one-sided orthogonalization strategy is more aggressive, as it orthogonalizes just one basis instead of two. That is why in our solver one-sided orthogonalization is an option that the user can activate, but by default full orthogonalization of the three bases is carried out. However, in our tests we have not found any situation where the one-sided variant leads to large residual norms or a bad level of orthogonality of the computed bases.

If we consider the overall computational cost, the saving from one-sided orthogonalization may be modest, since each iteration of the loop of \cref{alg:kilmer} performs a least squares solve at \cref{alg:kilmer:expand2}, whose cost is usually much higher, especially if using the LSQR iterative method. See results in \cref{sec:results}.

\section{Details of the implementation in SLEPc}\label{sec:implem}

We now discuss a few relevant details of our implementation, with which the results shown in \cref{sec:results} have been obtained. The implementation is included in one of the solvers of SLEPc, the Scalable Library for Eigenvalue Problem Computations~\cite{Hernandez:2005:SSF}. We first give an overview of this library, before going into the details.

SLEPc is a parallel library intended for solving large-scale eigenvalue problems, mainly by iterative methods, that is, in cases where only part of the spectrum is required. It consists of several modules, each one for a different type of problem, including linear eigenproblems, polynomial eigenproblems, general nonlinear eigenproblems, and singular-value-type decompositions. The latter module is called \texttt{SVD} and is the one relevant for this work.

SLEPc can be seen as an extension of PETSc, the Portable, Extensible Toolkit for Scientific Computation~\cite{Balay:2022:PUM2}. PETSc provides a lot of functionality related to data structures and solvers useful for large-scale scientific applications modeled by partial differential equations. Apart from the basic data objects for working with matrices and vectors, the only PETSc modules that are relevant for this work are those implementing iterative methods for linear systems of equations (\texttt{KSP}) and preconditioners (\texttt{PC}).

The model of parallelism of PETSc/SLEPc is message-passing (MPI), where every object is created with respect to an MPI communicator, implying that certain operations are carried out collectively by all processes belonging to that communicator. A second level of parallelism is also available to further accelerate the local computations either via CPU threads or by launching kernels to a GPU, but we do not consider them in this paper.

The solver modules are essentially a collection of solvers with a common user interface, where the user can choose the solver to be employed. This selection can be done at run time via command-line options, which confers a lot of flexibility by allowing also to specify solver parameters such as the dimension of the subspace, the tolerance, and many more. In the case of the \texttt{SVD} module, it contains a number of solvers, from which we highlight the following:
\begin{itemize}
  \item \texttt{cross}. This solver originally contained a gateway for computing the SVD via the linear eigenvalue problem module, in particular with the equivalent cross product matrix eigenproblem, \eqref{eq:eigleft} or \eqref{eq:eigright}. In this work, we have extended this for the GSVD via the equivalent eigenproblem~\eqref{eq:gsvdeigcross}.
  \item \texttt{cyclic}. Similarly to the previous one, this solver had code for the SVD formulated via the cyclic eigenproblem~\eqref{eq:cyclic}, and we have extended it for the GSVD using the equivalent cyclic eigenproblems~\eqref{eq:gsvdeigcyclic}.
  \item \texttt{trlanczos}. This is the thick-restart Lanczos solver. The implementation for the SVD is described in~\cite{Roman:2022:SUM2}. The extension to the GSVD, following the methodology detailed in \cref{sec:trlanczos}, constitutes the main contribution of this paper.
\end{itemize}

At the core of the thick-restart Lanczos solver is the lower-upper\footnote{We have also implemented the upper-upper variant (Zha's method), including thick restart, but we do not consider this here to keep the presentation short.} joint bidiagonalization of \cref{alg:kilmer}. One of the main operations in this algorithm is the orthogonalization of vectors, which in SLEPc is carried out with an iterated classical Gram-Schmidt variant that is both efficient in parallel and numerically stable~\cite{Hernandez:2007:PAE}. Another notable operation is the expansion of the $\widetilde{V}$ Lanczos basis, which implies the solution of the least squares problem~\eqref{eq:leastsquares}. This is done using PETSc's functionality, as described next.

PETSc allows combining iterative methods from \texttt{KSP} with preconditioners from \texttt{PC}, for instance \texttt{gmres} and \texttt{jacobi}. Direct linear solvers can be used with a special \texttt{KSP} that means precondition only, e.g., \texttt{preonly} and \texttt{lu}. There are a limited number of instances of \texttt{KSP} and \texttt{PC} that can handle rectangular matrices, so the possibilities for least squares problems are:
\begin{itemize}
  \item A purely iterative method, i.e., without preconditioning, via the \texttt{KSP} solver \texttt{lsqr} that implements the LSQR method~\cite{Paige:1982:LAS}, or the \texttt{cgls} method, which is equivalent to applying the conjugate gradients on the normal equations.
  \item A direct method, with \texttt{preonly} and \texttt{qr}. This requires installing PETSc together with SuiteSparse, which implements a sparse QR factorization~\cite{Davis:2011:SQR}.
  \item A combination of the two, i.e., \texttt{lsqr} and \texttt{qr}. Since \texttt{qr} is a direct method, the \texttt{lsqr} method will finish in just one iteration. Note that LSQR requires that the preconditioner is built for the normal equations matrix $Z^TZ$, which means that in the case of \texttt{qr} the preconditioner application will consist in a forward triangular solve followed by a backward triangular solve with $R$, the computed triangular factor. Even though the number of iterations is one, due to how the method is implemented the number of preconditioner applications is two, but still this is usually cheaper than \texttt{preonly} with \texttt{qr} because the latter needs applying the orthogonal factor of the QR decomposition, which is avoided with \texttt{lsqr}.
  \item The iterative method \texttt{lsqr} with a parallel preconditioner for the normal equations. One possible parallel preconditioner is block Jacobi (\texttt{bjacobi}), a simple domain-decomposition preconditioner where each process applies a local preconditioner built from the block-diagonal part of the distributed matrix. In our case, the local preconditioner can be \texttt{qr}. A more effective preconditioner is algebraic multigrid as implemented in Hypre~\cite{Falgout:2002:HLH} or an advanced domain decomposition scheme as implemented in HPDDM~\cite{Jolivet:2021:KPE}. The latter has been specifically adapted for the case of the normal equations following the method in~\cite{Daas:2022:RAD}.
\end{itemize}

We will consider \texttt{lsqr} with or without preconditioner. The former case requires that matrix $Z$ is built explicitly, by stacking the user matrices $A$ and $B$, while the unpreconditioned \texttt{lsqr} can operate with  $A$ and $B$ independently. Hence, we have added a user parameter \texttt{explicitmatrix} that can be turned on if necessary. This option applies also to the \texttt{cross} and \texttt{cyclic} solvers to explicitly build the matrices of~\eqref{eq:gsvdeigcross} and~\eqref{eq:gsvdeigcyclic}, which is necessary in the case that a direct linear solver is to be used in the solution of the eigenproblem, e.g., to factorize matrix $B^TB$.

Apart from the joint bidiagonalization of \cref{alg:kilmer}, the main ingredient of the solver is the thick restart discussed in \cref{sec:trlanczos}. In SLEPc, the small-sized dense projected problem is handled within an auxiliary object called \texttt{DS}, with operations to solve the projected problem, sort the computed solution according to a certain criterion, and truncate it to a smaller size, among other. For the GSVD, we have implemented the computation of the CS decomposition~\eqref{eq-J-CS-decomp} using LAPACK's subroutine \texttt{\_GGSVD3}, and the truncation operation that manages the arrowhead shapes in \cref{fig:thick-gsvd} using a compact storage.

A common parameter in all SLEPc solvers is \texttt{ncv}, the number of column vectors, that is, the maximum allowed dimension of the subspaces. In SVD-type computations, the user can also choose whether to compute the largest (default) or smallest (generalized) singular values and vectors. The user interface of the GSVD thick-restart solver includes three additional parameters: the scale factor, discussed in \cref{sec:scaling}, a flag to enable one-sided orthogonalization, discussed in \cref{sec:oneside}, and the restart parameter, indicating the proportion of basis vectors that must be kept after restart (the default is 50\%).

A final comment is about locking. The absolute value of the spikes of arrows in \cref{fig:thick-gsvd} are precisely the residual bounds used in the convergence criterion~\eqref{eq-rgsvd-bound}. As soon as the generalized singular quadruples converge, the leading values of these spikes become small (below the tolerance). Then one could set these values to zero explicitly, with the consequent decoupling from the rest of the bidiagonalization. This is called locking, a form of deflation in which the converged solutions are used in the orthogonalization but are excluded from the rest of operations. In our solver, locking is active by default, but the user can deactivate it so that the full bidiagonalization is updated at restart (including converged vectors). All results discussed in next section use locking.

\section{Computational results}\label{sec:results}

In this section we present results of some computational experiments to assess the performance of the solvers in terms of accuracy, convergence, and parallel scalability. The executions have been carried out on the Tirant III computer, which consists of 336 computing nodes, each of them with two Intel Xeon SandyBridge E5-2670 processors (16 cores each) running at 2.6 GHz with 32 GB of memory, connected with an Infiniband network. We allocated 4 MPI processes per node at most. The presented results correspond to SLEPc version 3.18, together with PETSc 3.18\footnote{With a patch required to improve the performance when using the \texttt{hypre} preconditioner with \texttt{lsqr}, that will be included in PETSc 3.19.}, SuiteSparse 5.13, hypre 2.25 and MUMPS 5.5. All software has been compiled with Intel C and Fortran compilers (version 18) and Intel MPI.

\begin{table}%
\caption{Description of the test problems used in the computational experiments: number of rows $m$ and columns $n$ of $A$, number of rows $p$ of $B$, number of requested generalized singular values \texttt{nsv}, whether largest or smallest generalized singular values are wanted, and the first computed value.}
\label{tab:gsvdproblems}
\begin{minipage}{\columnwidth}
\begin{center}
\begin{tabular}{lcccccc}
\hline
name & $m$ & $n$ & $p$ & \texttt{nsv} & \texttt{which} & first value\\
\hline
\textsf{diagonal}   & 500000 & 500000& 500000 & 20 & largest  & $\sigma_1=0.57735$ \\
\textsf{invinterp}  & 262144 & 262144& 524288 & 20 & largest  & $\sigma_1=25701.6$  \\
\textsf{hardesty2}  & 626256 & 303645& 303645 & 20 & largest  & $\sigma_1=53990.8$ \\
\textsf{3elt}       &   4720 &  4720 & 4721   &  5 & largest  & $\sigma_1=8727.4$  \\
\textsf{gemat12}    &   4929 &  4929 & 4930   &  5 & largest  & $\sigma_1=47738$   \\
\textsf{onetone1}   &  36057 & 36057 & 36058  &  5 & smallest & $\sigma_n=3.44\cdot 10^{-4}$ \\
\hline
\end{tabular}
\end{center}
\end{minipage}
\end{table}%

\Cref{tab:gsvdproblems} lists the problems used in the experiments, summarizing some properties and parameters. Here is a short description of the problems:
\begin{itemize}
  \item The \textsf{diagonal} problem is taken from~\cite[Example 1]{Zha:1996:CGS}, but with larger dimension. The problem matrices are computed as $A=CD$ and $B=SD$, where $C=\diag(c_i)$ with $c_i=(n-i+1)/2n$, $S=\sqrt{I-C^2}$, and $D=\diag(d_i)$ with $d_i=\lceil 4i/n\rceil + r_i$ where $r_i$ is a random number uniformly distributed on $[0,1]$.
  \item The \textsf{invinterp} case corresponds to the inverse interpolation problem of IR Tools~\cite[\S3.2]{Gazzola:2018:ITM}, intended to test iterative regularization methods for linear inverse problems. The GSVD can be used in this context. In this case, the $A$ matrix comes from the interpolation relations of $n$ points at random locations with respect to a 2D regular grid of side $\sqrt{n}$, and $B$ is the regularization matrix. In our case, for matrix $B$ we use the last case listed in~\cite[\S4.2]{Gazzola:2018:ITM} (2D Laplacian with zero derivative enforced on one of the boundaries) except that we do not compute the compact QR decomposition of this matrix.
  \item In the \textsf{3elt}, \textsf{gemat12} and \textsf{onetone1} problems, the $A$ matrix is the matrix with the same name taken from the SuiteSparse matrix collection~\cite{Davis:2011:UFS}, while the $B$ matrix is a bidiagonal matrix with one more row than columns and values $-1$ and $1$ on the subdiagonal and diagonal, respectively. The \textsf{hardesty2} matrix also belongs to the SuiteSparse matrix collection, but instead of using it with a bidiagonal matrix, we take the bottom square block as matrix $B$ and the remaining top rows as matrix $A$.
\end{itemize}

\begin{table}%
\caption{Computational results for the first three test cases, solved sequentially and in parallel (with 16 MPI processes), using a sparse QR factorization or the LSQR iterative method for the least squares solves, respectively. We show the scale factor $\gamma$, the number of Lanczos restarts, the total number of iterations of the linear solver, the execution time in seconds, and the maximum relative residual norm of all computed solutions.}
\label{tab:gsvdresults}
\begin{minipage}{\columnwidth}
\begin{center}
\begin{tabular}{lccccccc}
\hline
test problem& LS solver&processes & $\gamma$&restarts&linear its&time& residual \\
\hline
\textsf{diagonal}   &QR&1   & 1 &    763 &12479& 2847 &$6\cdot 10^{-9}$ \\
\textsf{invinterp}  &QR        &1  & 50   &69  &1232& 585  &$2\cdot 10^{-9}$ \\
\textsf{hardesty2}  &QR&1   & 1000 &3  &102&81  &$6\cdot 10^{-11}$ \\
\hline
\textsf{diagonal}   &LSQR  &16  & 1 &818  & 701867 &489  &$6\cdot 10^{-9}$ \\
\textsf{invinterp}  &LSQR-hypre&16 & 50  &  65& 41049& 407  &$1\cdot 10^{-7}$ \\
\textsf{hardesty2}  &LSQR-hypre&16 & 1000  &  4& 26813 & 218 &$3\cdot 10^{-10}$ \\
\hline
\end{tabular}
\end{center}
\end{minipage}
\end{table}%

All results in this section are obtained with the \texttt{explicitmatrix} flag set, except for the \textsf{diagonal}, \textsf{3elt}, \textsf{gemat12} and \textsf{onetone1} problems when using LSQR without a preconditioner. The stopping tolerance for LSQR is $10^{-10}$.
\Cref{tab:gsvdresults} illustrates the performance of our solver with the first three test cases. In sequential runs we use the sparse QR factorization to solve the least squares problems, and in parallel we employ LSQR, which may need a lot of iterations (the \textsf{diagonal} and \textsf{invinterp} problems can be solved without preconditioner, but in \textsf{hardesty2} a preconditioner is required). From the table, we can see that both \textsf{invinterp} and \textsf{hardesty2} need more time in parallel than sequentially, because LSQR is slow (this is also related to the scale factor as discussed below). Still, parallelism allows solving very large problems where computing a sequential sparse QR is not viable.

To assess the accuracy of the computed generalized singular quadruples, we use the residual norm~\eqref{eq:residual} relative to the norm of the stacked matrix~\eqref{eq-QR}, $\left\|r^\mathrm{GSVD}\right\|_2/\|Z\|_\infty$. The default tolerance for the thick-restart Lanczos solver is $10^{-8}$, and we see in \cref{tab:gsvdresults} that the relative residual norm is below the tolerance in all cases except for \textsf{invinterp} with LSQR. This can be fixed by requesting a more stringent tolerance for the LSQR stopping criterion, otherwise a bad accuracy of the inner iteration prevents attaining the requested accuracy in the outer one.

In the case of a sparse QR factorization, the value of the total number of iterations of the linear solver in \cref{tab:gsvdresults} can be interpreted as the number of least squares problems that need to be solved. To understand this value, we have to take into account that the default value of \texttt{ncv} (number of column vectors) is equal to $\max\{2\cdot\texttt{nsv},10\}$, that is, $40$ in the case of the test problems analyzed in \cref{tab:gsvdresults}. In addition to the least squares problems required during the Lanczos iteration, the postprocessing to obtain the final right singular vectors $g_i=R^{-1}w_i$ requires additional least squares solves.

\begin{figure}
    \centering
    \begin{tikzpicture}
       \begin{axis}[
           xbar stacked,
           enlargelimits=0.15,
           legend pos=outer north east,
           legend cell align={left},
           legend style={font=\footnotesize},
           width=12cm, height=3.5cm, enlarge y limits=0.5,
           symbolic y coords={restarted, non-restarted},
           ytick=data,
       ]
        \addplot coordinates {  
            (120,restarted) (936,non-restarted)
        };
        \addplot coordinates {  
            (422,restarted) (268,non-restarted)
        };
        \addplot coordinates {  
            (43,restarted) (62.8,non-restarted)
        };
        \legend{Orthogonaliz.,Least squares,Other ops.}
    \end{axis}
    \end{tikzpicture}
    \caption{\label{fig:non-restarted}Comparison of the execution time (in seconds) of the non-restarted and thick-restarted methods when solving the \textsf{invinterp} problem with a sparse QR for the least squares using 1 MPI process. Other operations include the initial computation of the QR factorization, the solution of the projected problem~\eqref{eq-J-CS-decomp} and the update of bases \eqref{eq-restart-a}-\eqref{eq-restart-b} required at restart (and at the end of the algorithm).}
\end{figure}
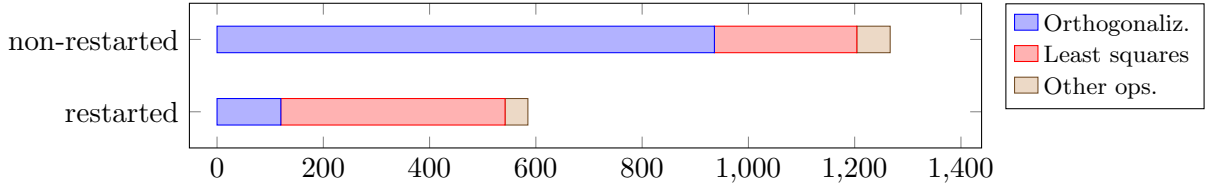

If we wanted to compare the performance of our thick-restarted GSVD solver (\cref{alg:restart}) with respect to a non-restarted solver (the algorithms in \cite{Zha:1996:CGS} and \cite{Jia:2021:JBP}, for instance), one possibility would be to emulate the latter by using the restarted solver with a sufficiently large basis size (\texttt{ncv}) so that restart is never required. We have carried out this comparison with the \textsf{invinterp} problem with a sparse QR for the least squares (second line of \cref{tab:gsvdresults}), using bases of \texttt{ncv}=765 vectors\footnote{This value has been chosen \emph{a priori}, but an actual implementation of a non-restarted solver should have to check convergence at every Lanczos step, which would increase the overall cost.}. Of course, in terms of memory use, the restarted solver is much cheaper since each of the three vector bases need to store only 40 vectors, compared to the 765 vectors of the non-restarted solver. But it is also cheaper in terms of computational cost, as shown in \cref{fig:non-restarted}. Clearly the orthogonalization cost blows up in the non-restarted variant. The cost associated with solving least squares problems increases in the thick-restart case, because the total number of performed least squares solves is 1232, larger than 789 required for the non-restarted version. Overall, the thick-restart solver takes less than half the time of the non-restarted one.

	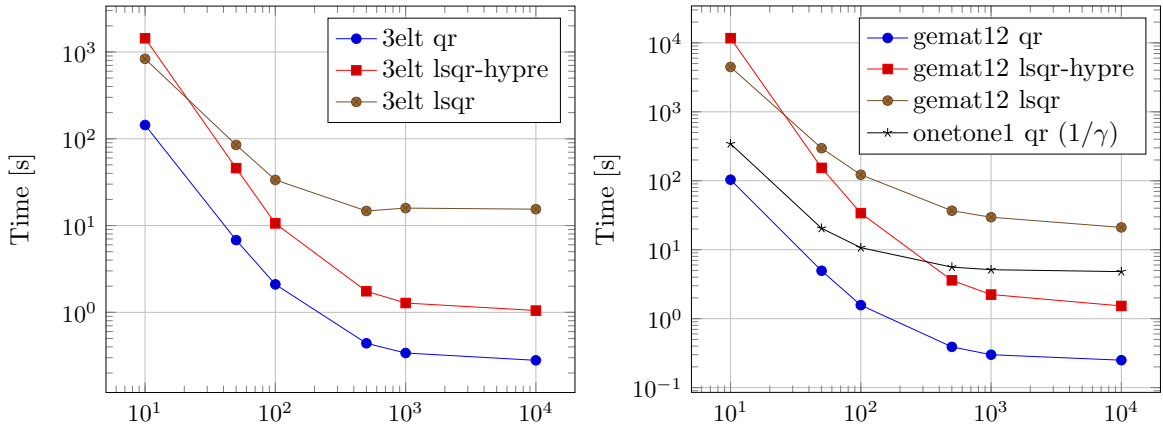
\begin{figure}
		\centering
		\begin{tikzpicture}[scale=0.9]
		\begin{loglogaxis}[
		ylabel={Time [s]},
		grid=major,
		ticklabel style={font=\small},
		legend style={anchor=north, cells={anchor=west}}, legend pos=north east,
		]
		\addplot coordinates { 
          (   10,  144.01)
          (   50,    6.80)
          (  100,    2.10)
          (  500,    0.44)
          ( 1000,    0.34)
          (10000,    0.28)
		};
		\addplot coordinates { 
          (   10, 1437.20)
          (   50,   45.87)
          (  100,   10.57)
          (  500,    1.75)
          ( 1000,    1.28)
          (10000,    1.05)
		};
		\addplot coordinates { 
          (   10,  833.85)
          (   50,   84.91)
          (  100,   33.55)
          (  500,   14.73)
          ( 1000,   15.89)
          (10000,   15.44)
		};
		\legend{3elt qr, 3elt lsqr-hypre, 3elt lsqr}
		\end{loglogaxis}
		\end{tikzpicture}
		\begin{tikzpicture}[scale=0.9]
		\begin{loglogaxis}[
		ylabel={Time [s]},
		grid=major,
		ticklabel style={font=\small},
		legend style={anchor=north, cells={anchor=west}}, legend pos=north east,
		]
		\addplot coordinates { 
          (   10,  103.58)
          (   50,    4.95)
          (  100,    1.57)
          (  500,    0.39)
          ( 1000,    0.30)
          (10000,    0.25)
		};
		\addplot coordinates { 
          (   10, 11689.00)
          (   50,  153.67)
          (  100,   33.89)
          (  500,    3.59)
          ( 1000,    2.24)
          (10000,    1.53)
		};
		\addplot coordinates { 
          (   10, 4466.70)
          (   50,  296.18)
          (  100,  122.20)
          (  500,   36.75)
          ( 1000,   29.57)
          (10000,   21.06)
		};
		\addplot [black,mark=star]coordinates { 
          (   10,  343.84)
          (   50,   20.43)
          (  100,   10.71)
          (  500,    5.59)
          ( 1000,    5.13)
          (10000,    4.81)
		};
		\legend{gemat12 qr,gemat12 lsqr-hypre, gemat12 lsqr, onetone1 qr ($1/\gamma$)}
		\end{loglogaxis}
		\end{tikzpicture}
		\caption{\label{fig:time-scale}Execution time (in seconds) with different scale factors, for matrices \textsf{3elt} (left), \textsf{gemat12} and \textsf{onetone1} (right).}
	\end{figure}

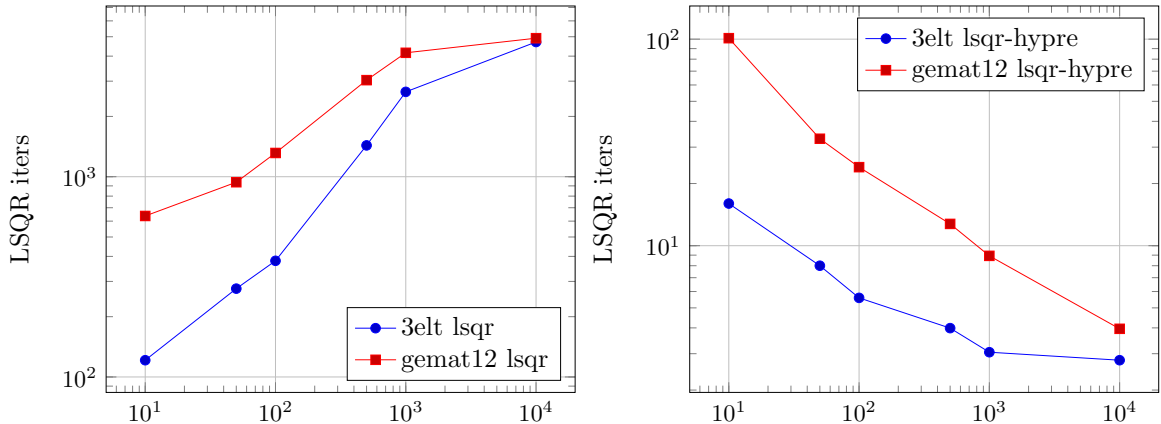
\begin{figure}
	\centering
		\begin{tikzpicture}[scale=0.9]
		\begin{loglogaxis}[
		ylabel={LSQR iters},
		grid=major,
		ticklabel style={font=\small},
		legend style={anchor=south east, cells={anchor=west}}, legend pos=south east
		]
		\addplot coordinates { %
          (   10,  121.23)
          (   50,  276.12)
          (  100,  380.22)
          (  500, 1433.43)
          ( 1000, 2652.77)
          (10000, 4706.33) 
		};
		\addplot coordinates { %
		  (   10,  636.76)
          (   50,  938.55)
          (  100, 1313.79)
          (  500, 3035.24)
          ( 1000, 4157.03)
          (10000, 4916.63)
		};
		\legend{3elt lsqr,gemat12 lsqr}
		\end{loglogaxis}
		\end{tikzpicture}
		\begin{tikzpicture}[scale=0.9]
		\begin{loglogaxis}[
		ylabel={LSQR iters},
		grid=major,
		ticklabel style={font=\small},
		legend style={anchor=south east, cells={anchor=west}}, legend pos=north east
		]
		\addplot coordinates { %
           (   10,   16.00)
           (   50,    8.00)
           (  100,    5.59)
           (  500,    3.99)
           ( 1000,    3.05)
           (10000,    2.79)
		};
		\addplot coordinates { %
           (   10,  101.09)
           (   50,   32.95)
           (  100,   24.00)
           (  500,   12.75)
           ( 1000,    8.93)
           (10000,    3.96)
		};
		\legend{3elt lsqr-hypre,gemat12 lsqr-hypre}
		\end{loglogaxis}
		\end{tikzpicture}
		\caption{\label{fig:lsqrits}Average number of LSQR iterations per least squares solve for \textsf{3elt} and \textsf{gemat12} with different scale factors, considering no preconditioner (left) or \texttt{hypre} preconditioner (right).}
	\end{figure}

The scale factor discussed in \cref{sec:scaling} may have a significant impact on the performance of the thick-restart Lanczos solver, as it has an influence on the number of Lanczos restarts as well as the number of iterations needed by LSQR. \Cref{fig:time-scale} compares the execution time of the last problems in \cref{tab:gsvdproblems} when the scale factor is changed, and using three different methods to solve the least squares problem: sparse QR, unpreconditioned LSQR and LSQR with an algebraic multigrid preconditioner from Hypre.
We can draw several conclusions.
First, in all these problems scaling is beneficial, and execution time is reduced whenever the scale factor $\gamma$ is increased, up to a point where it stabilizes. This is also true, but with the reciprocal of $\gamma$, if the smallest generalized singular values are computed (\texttt{onetone1} problem).
Second, in these problems it is much cheaper to use a sparse QR factorization for the least squares problems, compared to solving them with LSQR (either with or without preconditioning) as the latter may require many iterations (it does not converge in the \texttt{onetone1} problem). Still, LSQR is currently necessary for parallel runs, and in \cref{fig:lsqrits} we study how the number of iterations of LSQR changes with respect to the scale factor. As the scale factor increases, the unpreconditioned LSQR method has more difficulties to converge, but the preconditioned LSQR improves its convergence speed.

	\begin{figure}
	\centering
		\begin{tikzpicture}[scale=0.9]
		\begin{semilogyaxis}[
		log origin=infty,
		ylabel={Time [s]},
        ybar,
		bar width=6pt,
		xtick=data,
		xticklabels={invinterp,hardesty2,3elt,gemat12,onetone1},
		enlargelimits=0.2,
        major x tick style = transparent,
        x tick label style={rotate=45,anchor=east},
        legend style={anchor=south east, cells={anchor=west}}, legend pos=north east
		]
		\addplot coordinates { 
                  (1,326)
                  (2,81)
		  (3,0.28)
		  (4,0.25)
		  (5,4.81)
		};
		\addplot coordinates { 
                  (1,89)
                  (2,18)
		  (3,0.32)
		  (4,0.34)
		  (5,2.32)
		};
		\addplot coordinates { 
                  (1,183)
                  (2,728)
		  (3,1.45)
		  (4,1.37)
		  (5,11.06)
		};
		\legend{trlanczos,cross,cyclic}
		\end{semilogyaxis}
		\end{tikzpicture}
		\caption{\label{fig:solvers}Comparing the thick-restart Lanczos solver (using the best scale parameter $\gamma$ and sparse QR for the least squares problems) with the cross and cyclic solvers (using MUMPS for the linear solves).}
	\end{figure}
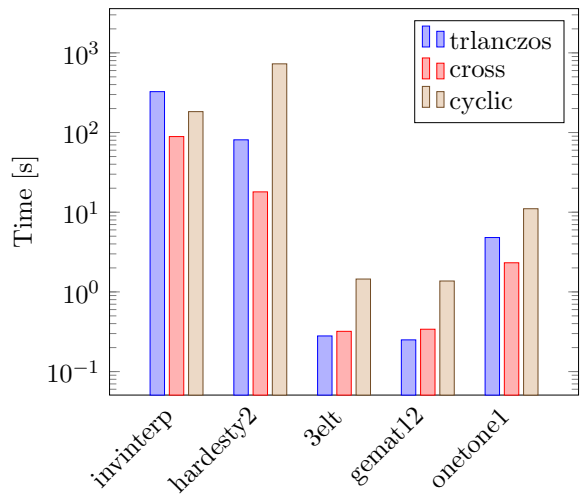

\Cref{fig:solvers} shows the execution time of the thick-restart Lanczos solver compared with the cross and cyclic solvers discussed in \cref{sec:equivep}. In all cases we use a direct method for the linear solves (SuiteSparse for Lanczos and MUMPS for cross and cyclic). There is no clear winner, but take into account that in Lanczos we are using the scale factor $\gamma$ that gave the smallest time in our tests for each problem, otherwise the Lanczos solver is not competitive in general.
However, it is worth noting that for the \textsf{hardesty2} problem, the Lanczos solver is the only one that provides a suitable solution in terms of the residual. In this case, the matrix $B^TB$ in the cross and cyclic methods is singular and, because the generalized symmetric-definite Lanczos eigensolver fails, one has to solve the equivalent eigenproblem as non-symmetric, resulting in much lower accuracy ($\approx1\cdot10^{-4}$), compared to the Lanczos GSVD solver ($\approx6\cdot10^{-11}$). In the other cases, the relative residual norm of the computed solutions is similar in all solvers.

	\begin{figure}
		\centering
\resizebox{\textwidth}{!}{
		\begin{tikzpicture}
		\begin{loglogaxis}[
		title={\textsf{diagonal} - total time},
		ylabel={Time [s]},
		grid=major,
		log basis x=2,
		xtick={1,2,4,8,16,32,64},
		xticklabels={1,2,4,8,16,32,64},
		ticklabel style={font=\small},
legend style={draw=none, legend columns=4, legend to name=legendname2,font=\small,cells={anchor=west}}
		]
		\addplot coordinates { %
			(  1, 7.4507e+03)
			(  2, 3.6038e+03)
			(  4, 1.9231e+03)
			(  8, 8.0392e+02)
			( 16, 4.1416e+02)
			( 32, 2.1563e+02)
			( 64, 1.3569e+02)
		};
		\addplot coordinates { %
			(  1, 6.6613e+03)
			(  2, 3.3257e+03)
			(  4, 1.6989e+03)
			(  8, 7.4470e+02)
			( 16, 3.6732e+02)
			( 32, 1.9930e+02)
			( 64, 1.1667e+02)
		};
    \addplot[color=gray, very thick,dashed, mark options={solid}] coordinates {
        (2, 1000) (4, 500) (8, 250) (16, 125)
    };
		\legend{twoside,oneside,ideal scaling}
		\end{loglogaxis}
		\end{tikzpicture}
\hfill
		\begin{tikzpicture}
		\begin{loglogaxis}[
		title={\textsf{diagonal} - orthogonalization only},
		grid=major,
		log basis x=2,
		xtick={1,2,4,8,16,32,64},
		xticklabels={1,2,4,8,16,32,64},
		ticklabel style={font=\small},
		]
		\addplot coordinates { %
			(  1, 1.3548e+03)
			(  2, 6.4835e+02)
			(  4, 3.5231e+02)
			(  8, 1.6582e+02)
			( 16, 7.6805e+01)
			( 32, 3.6858e+01)
			( 64, 1.8521e+01)
		};
		\addplot coordinates { %
			(  1, 5.2097e+02)
			(  2, 2.6152e+02)
			(  4, 1.3400e+02)
			(  8, 7.0736e+01)
			( 16, 2.7122e+01)
			( 32, 1.2542e+01)
			( 64, 6.0540e+00)
		};
    \addplot[color=gray, very thick,dashed, mark options={solid}] coordinates {
        (2, 100) (4, 50) (8, 25) (16, 12.5)
    };
		\end{loglogaxis}
		\end{tikzpicture}
}
\ref{legendname2}
		\caption{\label{fig:paraldiag}Execution time (in seconds) with up to 64 MPI processes for the  \textsf{diagonal} problem, solved with thick-restart Lanczos using one-sided or two-sided orthogonalization. The left plot shows the total time, while the right one accounts only for the time of orthogonalization.}
	\end{figure}
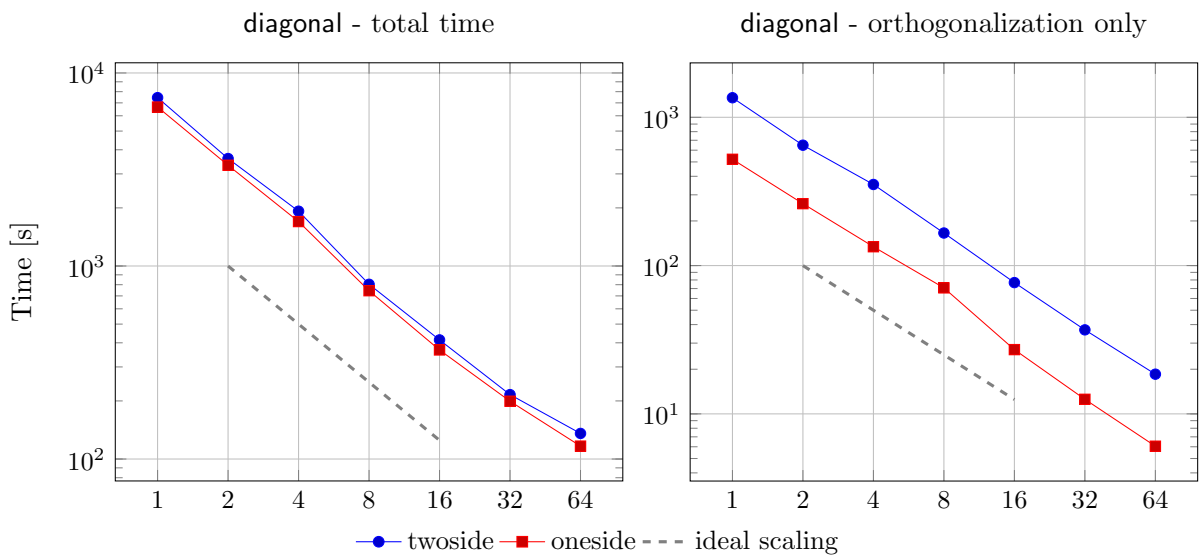

	\begin{figure}
		\centering
		\begin{tikzpicture}[scale=0.9]
		\begin{loglogaxis}[
		ylabel={Time [s]},
		grid=major,
		log basis x=2,
		xtick={1,2,4,8,16,32,64},
		xticklabels={1,2,4,8,16,32,64},
		ticklabel style={font=\small},
		]
		\addplot coordinates { %
			(  1, 4.7584e+03)
			(  2, 2.0769e+03)
			(  4, 1.2205e+03)
			(  8, 6.7490e+02)
			( 16, 3.9775e+02)
			( 32, 2.8252e+02)
			( 64, 2.8718e+02)
		};
		\end{loglogaxis}
		\end{tikzpicture}
		\caption{\label{fig:irtools}Execution time (in seconds) with up to 64 MPI processes for the \texttt{invinterp} problem, solved with thick-restart Lanczos and \texttt{hypre}-preconditioned LSQR.}
	\end{figure}
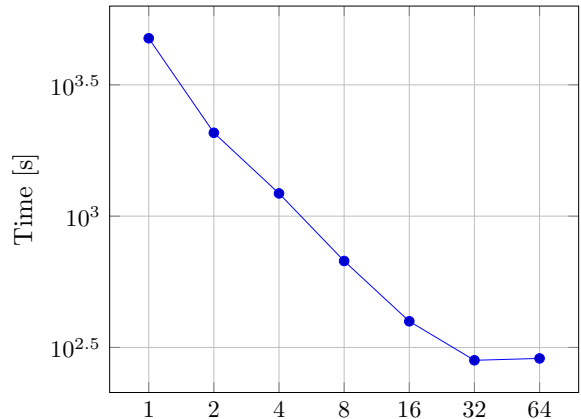

We conclude this section by analyzing parallel performance. \Cref{fig:paraldiag} plots the parallel execution time of the Lanczos solver for the \textsf{diagonal} problem with up to 64 MPI processes, using unpreconditioned LSQR. We can see that the scaling is very good, close to the ideal one. The figure shows the one-sided variant together with the default one (two-sided). The one-sided variant is always faster, but the difference is not too significant because orthogonalization amounts to only a modest percentage of the overall cost. To better appreciate the gain of one-sided orthogonalization, the right panel of \cref{fig:paraldiag} shows only the time of the orthogonalization, with a factor of about $2.5$ improvement of the one-sided scheme with respect to the default.

The \textsf{diagonal} problem is very favorable for parallel computing, as the matrix-vector product is trivially parallelizable. On the other extreme, the \textsf{invinterp} problem has a very disadvantageous sparsity pattern (nonzeros located at random positions), so the speedup shown in \cref{fig:irtools} is good up to 32 MPI processes, but stagnates afterwards.

\section{Concluding remarks}\label{sec:concl}

We have developed a thick restart mechanism for the joint Lanczos bidiagonalization to compute a partial GSVD of a large-scale matrix pair. This is a very important piece to make the solver usable in the context of real applications. We have developed a fully-fledged implementation in the SLEPc library, that in addition to the restart, includes other interesting features such as one-sided orthogonalization, scaling, or locking. The solver is very flexible, in the sense that the user can indicate at run time whether the associated least squares problems must be solved with a direct or iterative method, as well as specify many other settings such as the dimension of the Krylov subspace or the scale factor.

We have conducted a number of computational experiments, showing that our solver is numerically robust, computationally efficient and scalable in parallel runs. If an appropriate scale factor is used, the performance of Lanczos method is on a par with that of the cross and cyclic solvers for the GSVD, which we have also developed during this work.


\begin{acknowledgements}
The authors would like to thank P.\ Jolivet for useful suggestions about parallel preconditioning in PETSc.
The computational experiments of \cref{sec:results} were carried out on the supercomputer Tirant III belonging to Universitat de Val\`encia.
\end{acknowledgements}

\bibliographystyle{plain}

\begin{thebibliography}{10}

\bibitem{Anderson:1999:LUG}
E.~Anderson, Z.~Bai, C.~Bischof, L.~S. Blackford, J.~Demmel, J.~Dongarra, J.~Du
  Croz, A.~Greenbaum, S.~Hammarling, A.~McKenney, and D.~Sorensen.
\newblock {\em {LAPACK} Users' Guide}.
\newblock Society for Industrial and Applied Mathematics, Philadelphia, PA,
  third edition, 1999.

\bibitem{Baglama:2005:AIR}
James Baglama and Lothar Reichel.
\newblock Augmented implicitly restarted {Lanczos} bidiagonalization methods.
\newblock {\em {SIAM} J. Sci. Comput.}, 27(1):19--42, 2005.

\bibitem{Bai:2000:TSA}
Z.~Bai, J.~Demmel, J.~Dongarra, A.~Ruhe, and H.~van~der Vorst, editors.
\newblock {\em Templates for the Solution of Algebraic Eigenvalue Problems: A
  Practical Guide}.
\newblock Society for Industrial and Applied Mathematics, Philadelphia, PA,
  2000.

\bibitem{Balay:2022:PUM2}
S.~Balay, S.~Abhyankar, M.~F. Adams, J.~Brown S.~Benson, P.~Brune,
  K.~Buschelman, E.~M. Constantinescu, L.~Dalcin, A.~Dener, V.~Eijkhout,
  J.~Faibussowitsch, W.~D. Gropp, V.~Hapla, T.~Isaac, P.~Jolivet, D.~Karpeev,
  D.~Kaushik, M.~G. Knepley, F.~Kong, S.~Kruger, D.~A. May, L.~Curfman McInnes,
  R.~Tran Mills, L.~Mitchell, T.~Munson, J.~E. Roman, K.~Rupp, P.~Sanan,
  J.~Sarich, B.~F. Smith, S.~Zampini, H.~Zhang, H.~Zhang, and J.~Zhang.
\newblock {\em {PETSc/TAO} Users Manual Revision 3.18}, 2022.
\newblock Argonne Technical Memorandum ANL-21/39.

\bibitem{Campos:2016:RQM}
C.~Campos and J.~E. Roman.
\newblock Restarted {Q-Arnoldi-type} methods exploiting symmetry in quadratic
  eigenvalue problems.
\newblock {\em {BIT}}, 56(4):1213--1236, 2016.

\bibitem{Daas:2022:RAD}
H.~Al Daas, P.~Jolivet, and J.~A. Scott.
\newblock A robust algebraic domain decomposition preconditioner for sparse
  normal equations.
\newblock {\em {SIAM} J. Sci. Comput.}, 44(3):A1047--A1068, 2022.

\bibitem{Davis:2011:SQR}
T.~A. Davis.
\newblock Algorithm 915, {SuiteSparseQR}: Multifrontal multithreaded
  rank-revealing sparse {QR} factorization.
\newblock {\em {ACM} Trans. Math. Software}, 38(1):1--22, 2011.

\bibitem{Davis:2011:UFS}
T.~A. Davis and Y.~Hu.
\newblock The {University of Florida Sparse Matrix Collection}.
\newblock {\em {ACM} Trans. Math. Software}, 38(1):1:1--1:25, 2011.

\bibitem{Edelman:2020:GWA}
A.~Edelman and Y.~Wang.
\newblock The {GSVD}: Where are the ellipses?, matrix trigonometry, and more.
\newblock {\em {SIAM} J. Matrix Anal. Appl.}, 41(4):1826--1856, 2020.

\bibitem{Falgout:2002:HLH}
Robert~D. Falgout and Ulrike~Meier Yang.
\newblock hypre: {A} library of high performance preconditioners.
\newblock In Peter M.~A. Sloot, Chih Jeng~Kenneth Tan, Jack Dongarra, and
  Alfons~G. Hoekstra, editors, {\em Computational Science - {ICCS} 2002,
  International Conference, Amsterdam, The Netherlands, April 21-24, 2002.
  Proceedings, Part {III}}, volume 2331 of {\em Lect. Notes Comp. Sci.}, pages
  632--641. Springer, 2002.

\bibitem{Gazzola:2018:ITM}
S.~Gazzola, P.~C. Hansen, and J.~G. Nagy.
\newblock {IR Tools}: a {MATLAB} package of iterative regularization methods
  and large-scale test problems.
\newblock {\em Numer. Algorithms}, 81(3):773--811, 2018.

\bibitem{Golub:1996:MC}
G.~H. Golub and C.~F. van Loan.
\newblock {\em Matrix Computations}.
\newblock The Johns Hopkins University Press, Baltimore, MD, third edition,
  1996.

\bibitem{Hernandez:2007:PAE}
V.~Hernandez, J.~E. Roman, and A.~Tomas.
\newblock Parallel {Arnoldi} eigensolvers with enhanced scalability via global
  communications rearrangement.
\newblock {\em Parallel Comput.}, 33(7--8):521--540, 2007.

\bibitem{Hernandez:2005:SSF}
V.~Hernandez, J.~E. Roman, and V.~Vidal.
\newblock {SLEPc}: A scalable and flexible toolkit for the solution of
  eigenvalue problems.
\newblock {\em {ACM} Trans. Math. Software}, 31(3):351--362, 2005.

\bibitem{Hernandez:2008:REP}
Vicente Hernandez, Jose~E. Roman, and Andres Tomas.
\newblock A robust and efficient parallel {SVD} solver based on restarted
  {Lanczos} bidiagonalization.
\newblock {\em Electron. Trans. Numer. Anal.}, 31:68--85, 2008.

\bibitem{Hochstenbach:2009:JTM}
M.~E. Hochstenbach.
\newblock A {Jacobi--Davidson} type method for the generalized singular value
  problem.
\newblock {\em Linear Algebra Appl.}, 431(3-4):471--487, 2009.

\bibitem{Huang:2020:CFC}
J.~Huang and Z.~Jia.
\newblock On choices of formulations of computing the generalized singular
  value decomposition of a large matrix pair.
\newblock {\em Numer. Algorithms}, 87(2):689--718, 2020.

\bibitem{Jia:2021:JBP}
Z.~Jia and H.~Li.
\newblock The joint bidiagonalization process with partial reorthogonalization.
\newblock {\em Numer. Algorithms}, 88(2):965--992, 2021.

\bibitem{Jia:2023:JBM}
Z.~Jia and H.~Li.
\newblock The joint bidiagonalization method for large {GSVD} computations in
  finite precision.
\newblock {\em {SIAM} J. Matrix Anal. Appl.}, 44(1):382--407, 2023.

\bibitem{Jia:2020:JBB}
Z.~Jia and Y.~Yang.
\newblock A joint bidiagonalization based iterative algorithm for large scale
  general-form {Tikhonov} regularization.
\newblock {\em App. Numer. Math.}, 157:159--177, 2020.

\bibitem{Jolivet:2021:KPE}
P.~Jolivet, J.~E. Roman, and S.~Zampini.
\newblock {KSPHPDDM} and {PCHPDDM}: Extending {PETSc} with advanced {Krylov}
  methods and robust multilevel overlapping {Schwarz} preconditioners.
\newblock {\em Comput. Math. Appl.}, 84:277--295, 2021.

\bibitem{Kilmer:2007:PAG}
M.~E. Kilmer, P.~C. Hansen, and M.~I. Espa{\~{n}}ol.
\newblock A projection-based approach to general-form {Tikhonov}
  regularization.
\newblock {\em {SIAM} J. Sci. Comput.}, 29(1):315--330, 2007.

\bibitem{Paige:1981:TGS}
C.~C. Paige and M.~A. Saunders.
\newblock Towards a generalized singular value decomposition.
\newblock {\em {SIAM} J. Numer. Anal.}, 18(3):398--405, 1981.

\bibitem{Paige:1982:LAS}
C.~C. Paige and M.~A. Saunders.
\newblock {LSQR}: an algorithm for sparse linear equations and sparse least
  squares.
\newblock {\em {ACM} Trans. Math. Software}, 8(1):43--71, 1982.

\bibitem{Parlett:1980:SEP}
B.~N. Parlett.
\newblock {\em The Symmetric Eigenvalue Problem}.
\newblock Prentice-Hall, Englewood Cliffs, NJ, 1980.
\newblock Reissued with revisions by SIAM, Philadelphia, 1998.

\bibitem{Roman:2022:SUM2}
J.~E. Roman, C.~Campos, L.~Dalcin, E.~Romero, and A.~Tomas.
\newblock {SLEPc} users manual.
\newblock Technical Report DSIC-II/24/02--Revision 3.18, D. Sistemes
  Inform\`atics i Computaci\'o, Universitat Polit\`ecnica de Val\`encia, 2022.

\bibitem{Simon:2000:LMA}
Horst~D. Simon and Hongyuan Zha.
\newblock Low-rank matrix approximation using the {Lanczos} bidiagonalization
  process with applications.
\newblock {\em {SIAM} J. Sci. Comput.}, 21(6):2257--2274, 2000.

\bibitem{Loan:1976:GSV}
C.~F. van Loan.
\newblock Generalizing the singular value decomposition.
\newblock {\em {SIAM} J. Numer. Anal.}, 13(1):76--83, 1976.

\bibitem{Wu:2000:TLM}
Kesheng Wu and Horst Simon.
\newblock Thick-restart {Lanczos} method for large symmetric eigenvalue
  problems.
\newblock {\em {SIAM} J. Matrix Anal. Appl.}, 22(2):602--616, 2000.

\bibitem{Zha:1996:CGS}
H.~Zha.
\newblock Computing the generalized singular values/vectors of large sparse or
  structured matrix pairs.
\newblock {\em Numer. Math.}, 72(3):391--417, 1996.

\end{thebibliography}

\end{document}